\documentclass{elsarticle}

\usepackage{fullpage}
\usepackage{graphicx}
\usepackage{amsmath}
\usepackage{amsthm}
\usepackage{amssymb}
\usepackage{caption}
\usepackage{bm}
\usepackage[ruled,noend,norelsize]{algorithm2e}
\usepackage{natbib}
\usepackage{color}
\usepackage{verbatim}
\usepackage{booktabs}
\usepackage{subcaption}
\usepackage[bookmarks=false,colorlinks]{hyperref}
\AtBeginDocument{%
\hypersetup{
linkcolor=blue,   
citecolor=red,
urlcolor=magenta}}

\newcommand{\blue}{\textsc{Blue}\xspace}
\newcommand{\green}{\textsc{Green}\xspace}
\newcommand{\conv}{\textsc{Conv}\xspace}
\newcommand{\X}{\mathcal{X}}
\newcommand{\N}{\mathcal{N}}
\newcommand{\G}{\mathcal{G}}
\newcommand{\A}{\mathcal{A}}
\newcommand{\Ao}{\A_o}
\newcommand{\As}{\A_s}
\newcommand{\Ar}{\A_r}

\newcommand{\Y}{\mathcal{Y}}
\renewcommand{\P}{\mathcal{P}}
\newcommand{\D}{\mathcal{D}}
\newcommand{\R}{\mathcal{R}}
\newcommand{\V}{\mathcal{V}}
\newcommand{\C}{\mathcal{C}}
\newcommand{\M}{\mathcal{M}}
\newcommand{\AM}{\mathbb{M}}
\newcommand{\inc}{\gamma^{-}}
\newcommand{\out}{\gamma^{+}}

\newcommand{\dvvp}{\delta_{vv'}}
\newcommand{\dvpv}{\delta_{v'v}}
\newcommand{\minU}{\underline{U}}
\newcommand{\maxU}{\overline{U}}

\newcommand{\Vt}{\V^n(t)}
\newcommand{\Vit}{\V^n_{i}(t)}
\newcommand{\Vijt}{\V^n_{ij}(t)}
\renewcommand{\Re}{\mathbb{R}}

\renewcommand{\emptyset}{\varnothing}
\newcommand{\p}{\rho}
\newcommand{\usi}{\bar{s}_i}

\newcommand{\usij}{\bar{s}_{ij}}
\newcommand{\dt}{\Delta t}
\newcommand{\xt}{\bm{x}(t)}
\newcommand{\xtt}{\bm{x}(t+\dt)}
\newcommand{\xit}{x_i(t)}
\newcommand{\xjt}{x_j(t)}

\newcommand{\xitt}{x_i(t+\dt)}
\newcommand{\xjtt}{x_j(t+\dt)}
\newcommand{\xijt}{x_{ij}(t)}

\newcommand{\yit}{y_i(t)}
\newcommand{\yits}{y^\star_i(t)}

\newcommand{\yt}{\bm{y}(t)}
\newcommand{\Yt}{\bm{Y}(t)}
\newcommand{\Yit}{Y_i(t)}
\newcommand{\Yjt}{Y_j(t)}
\newcommand{\yijt}{y_{ij}(t)}
\newcommand{\at}{\bm{\alpha}(t)}
\newcommand{\atn}{\bm{\alpha}^n(t)}

\newcommand{\aijt}{\alpha_{ij}(t)}
\newcommand{\bt}{\bm{\beta}(t)}
\newcommand{\bijt}{\beta_{ij}(t)}
\newcommand{\bijpt}{\beta_{i'j'}(t)}
\newcommand{\pijt}{P_{ij}(t)}
\newcommand{\pij}{p_{ij}}

\newcommand{\pit}{\phi_i(t)}

\newcommand{\wit}{w_i(t)}

\newcommand{\mij}{\mu_{ij}(t)}
\newcommand{\mijp}{\mu_{i'j'}(t)}
\newcommand{\chit}{\chi_i(t)}
\newcommand{\chijt}{\chi_{ij}(t)}
\newcommand{\lijt}{\lambda_{ij}(t)}
\newcommand{\lijpt}{\lambda_{ij}^{i'j'}(t)}

\newcommand{\fii}{f_i}

\newcommand{\di}{d_i}
\renewcommand{\dj}{d_j}

\newcommand{\dit}{D_i(t)}
\newcommand{\djt}{D_j(t)}
\newcommand{\dxt}{\Delta \xt}
\newcommand{\dxit}{\Delta \xit}
\newcommand{\dxjt}{\Delta \xjt}
\newcommand{\zblue}{Z^n_B(\xt)}
\newcommand{\zgreen}{Z^n_G(\xt)}
\newcommand{\atng}{\bm{\alpha}_G^n(t)}
\newcommand{\atnb}{\bm{\alpha}_B^n(t)}
\newcommand{\btng}{\bm{\beta}_G^n(t)}
\newcommand{\btnb}{\bm{\beta}_B^n(t)}
\newcommand{\cijij}{c_{ij}^{i'j'}}
\newcommand{\gijij}{g_{ij}^{i'j'}}

\newcommand{\gt}{\bm{\gamma}(t)}

\newcommand{\git}{\gamma_{i}(t)}

\newcommand{\ghi}{\hat{\gamma}_i}
\newcommand{\gh}{\hat{\bm{\gamma}}}

\newcommand{\pipt}{\phi_{i'}(t)}
\newcommand{\xipt}{x_{i'}(t)}
\newcommand{\pijp}{p_{i'j'}}
\newcommand{\usijp}{\bar{s}_{i'j'}}

\DeclareMathOperator*{\argmax}{arg\,max}

\newtheorem{theorem}{Theorem}

\newtheorem{corollary}{Corollary}

\newtheorem{defi}{Definition}

\begin{document}
\allowdisplaybreaks[1]

\begin{frontmatter}

\title{Blue Phase: Optimal Network Traffic Control for Legacy and Autonomous Vehicles}
\author{David Rey and Michael W. Levin}

\begin{abstract}
With the forecasted emergence of autonomous vehicles in urban traffic networks, new control policies are needed to leverage their potential for reducing congestion. While several efforts have studied the fully autonomous traffic control problem, there is a lack of models addressing the more imminent transitional stage wherein legacy and autonomous vehicles share the urban infrastructure. We address this gap by introducing a new policy for stochastic network traffic control involving both classes of vehicles. We conjecture that network links will have dedicated lanes for autonomous vehicles which provide access to traffic intersections and combine traditional green signal phases with autonomous vehicle-restricted signal phases named \emph{blue} phases. We propose a new pressure-based, decentralized, hybrid network control policy that activates selected movements at intersections based on the solution of mixed-integer linear programs. We prove that the proposed policy is stable, \emph{i.e.} maximizes network throughput, under conventional travel demand conditions. We conduct numerical experiments to test the proposed policy under varying proportions of autonomous vehicles. Our experiments reveal that considerable trade-offs exist in terms of vehicle-class travel time based on the level of market penetration of autonomous vehicles. Further, we find that the proposed hybrid network control policy improves on traditional green phase traffic signal control for high levels of congestion, thus helping in quantifying the potential benefits of autonomous vehicles in urban networks.
\end{abstract}

\end{frontmatter}

%%%%%%%%%%%%%%%%%%%%%%%%%%%%%%%%%%%%%%%%%%%%%%%%%%%%%%%%%%%%%%%%%%%%%%%%%%%%
\section{Introduction}
%%%%%%%%%%%%%%%%%%%%%%%%%%%%%%%%%%%%%%%%%%%%%%%%%%%%%%%%%%%%%%%%%%%%%%%%%%%%

The emergence of autonomous vehicles (AVs) in urban networks leads to new operational challenges that have received a growing attention over the past few years. Of particular interest is the future paradigm wherein legacy (or human-operated) vehicles have been fully replaced with AVs. In this futuristic context, several studies have shown that the management of urban networks may benefit from forecasted technological advancements, such as improved trajectory control, automatic collision avoidance or platooning. However, the intermediate traffic state wherein legacy and autonomous vehicles co-exist has not been nearly as much examined by researchers. In such an intermediate traffic context, urban networks will need to adapt to make the most out of AV technology and allow the transition towards fully autonomous traffic, if this is ever to happen. The first AVs evolving in urban traffic networks are likely to have to abide by the existing infrastructure and legislation. However, this picture may change rapidly. We conjecture that with the increase in AV demand, urban traffic networks will adapt and that AV-specific infrastructure will be available to improve network operations at traffic intersections. 

In this paper, we hypothesize a hybrid traffic context wherein legacy vehicles (LVs) and AVs both have access to dedicated infrastructure in an urban transport network. We propose a new stochastic hybrid network control model to manage traffic intersections by combining traditional green signal phases with AV-restricted signal phases called \emph{blue} phases. During blue phases, only AVs are allowed to proceed through the intersection. On the other hand, green phases can be used by any type of vehicle. To avoid first-in-first-out restrictions on queue departures, we assume that AVs have dedicated lanes to access traffic intersections, hereby referred to as AV-lanes. Legacy vehicle lanes (LV-lanes) can be used by both LVs and AVs to access traffic intersection but for convenience we assume that AVs choose AV-lanes by default. We assume fixed route choice or known turning ratios and we prove that the proposed network control policy is stable under conventional travel demand conditions, \emph{i.e.} that queues are bounded. 

We first discuss the literature on network traffic control in transportation before investigating the state of the art in intersection management with AVs. We then position our paper with respect to the field and highlight our contributions.

%%%%%%%%%%%%%%%%%%%%%%%%%%%%%%%%%%%%%%%%%%%%%%%%%%%%%%%%%%%%%%%%%%%%%%%%%%%%
\subsection{Network Traffic Control}
%%%%%%%%%%%%%%%%%%%%%%%%%%%%%%%%%%%%%%%%%%%%%%%%%%%%%%%%%%%%%%%%%%%%%%%%%%%%

Traffic signal timing has been studied for decades, and some of the classic methods such as \citet{webster1958traffic}'s delay formula are widely known today. Optimization of individual signals is well-established, but optimal coordination of signals across networks is more complex~\citep{gartner1975optimization,robertson1991optimizing} as complete two-way coordination even through grid networks can be impossible to achieve. Signal timing is further complicated by additions such as transit signal priority~\citep{skabardonis2000control,dion2002rule}. Nevertheless, several commercial software tools (\emph{e.g.} Synchro) are currently available to assist city engineers with signal timings.

The literature on network traffic control spans across several fields such as optimization and control theory, transportation engineering and computer science. The seminal work of \citet{tassiulas1992stability} pioneered the research on stability conditions in network traffic control with an application to data packets routing in communication networks. The authors notably introduced the concept of back-pressure algorithm as method for decentralized control of a network of routers. In the context of urban transport networks, there is an extensive body of research on network traffic signal control. \citet{varaiya2013max} proposed a network traffic control policy based on max-pressure, a variant of the back-pressure algorithm wherein the control policy chooses the signal phase that maximizes the pressure at each intersection. Stability was proven assuming each turning movement has a distinct queue. \citet{wongpiromsarn2012distributed} proposed a pressure-based policy which maximizes throughput for unbounded queues. However, practical limitations such as link length require a careful choice of the pressure function to avoid queue spillback. Building on this effort, \citet{xiao2014pressure} proposed a pressure-releasing policy that accounts for finite queue capacities. Nonetheless, to more canonically apply the pressure-based routing they assumed that each turning movement has a separate queue, which is often not realistic. \citet{le2015decentralized} proposed a fixed-time policy wherein all phases are assigned a non-zero activation time and proved stability for fixed turning proportions and unbounded queues. Recently, \citet{valls2016convex} propose a convex optimization approach to traffic signal control that decouples the stability of the system from the choice of traffic control policy.

Starting from the seminal work of \citet{smith1979traffic}, several efforts have also attempted to model the impact of traffic signal optimization on route choice \citep{zhang2012traffic, gregoire2013back,zaidi2016back}. Other studies used bi-level optimization for signal timings that assume a user equilibrium behavior~\cite{yang1995traffic, sun2006bi}. \citet{le2017utility} used utility functions in max-pressure control to influence routing. However, while \citet{tassiulas1992stability}'s policy is provably throughput-optimal, max-pressure route choice makes no guarantees on the efficiency of the travel times. We assume fixed route choice in our network traffic control model to focus on green and blue phases signal control, and leave route choice for later studies.

%%%%%%%%%%%%%%%%%%%%%%%%%%%%%%%%%%%%%%%%%%%%%%%%%%%%%%%%%%%%%%%%%%%%%%%%%%%%
\subsection{Intersection Management with Autonomous Vehicles}
%%%%%%%%%%%%%%%%%%%%%%%%%%%%%%%%%%%%%%%%%%%%%%%%%%%%%%%%%%%%%%%%%%%%%%%%%%%%

With the advent of connected and autonomous vehicular technology over the past decades, researchers have progressively explored and proposed new intersection management paradigms instead of the traditional tricolor signal control system. \citet{dresner2004multiagent,dresner2005multiagent} proposed the autonomous intersection manager (AIM) protocol as an alternative to traffic signals for AVs. In AIM, vehicles use wireless communication channels to request a reservation from the intersection manager (IM). A reservation specifies the turning movement as well as the time at which the vehicle can enter the intersection. The IM simulates reservation requests on a space-time grid of tiles and accepts requests if they do not conflict with other reservations. If two vehicles will occupy the same tile at the same time, a potential conflict exists and the request must be rejected. The AIM protocol can be used with a First-Come-First-Served (FCFS) policy, wherein vehicles are prioritized based on their arrival time at the intersection. \citet{fajardo2011automated} and \citet{li2013modeling} compared FCFS with optimized traffic signals, and found that FCFS reduced delays in certain scenarios. However, \citet{levin2015paradoxes} found that FCFS created more delay with signals in certain common cases (\emph{e.g.} asymmetric intersections). Traffic signals are within the feasible region of AIM controls~\citep{dresner2007sharing}, so AIM can always perform at least as well as signals when optimized.

Naturally, FCFS is only one of many potential policies for reservations. \citet{schepperle2007agent,schepperle2008traffic} allowed vehicles to bid for priority: in addition to reservation requests, vehicles could also communicate their willingness-to-pay and receive priority access accordingly. The authors found that auctions reduced delay weighted by vehicle value-of-time, and \citet{vasirani2010market} found similar results for a network of intersections. \citet{carlino2013auction} used system bids to further improve the efficiency of auctions. However, \citet{levin2015intersection} found that high value-of-time vehicles could become trapped behind low value-of-time vehicles, making it difficult to achieve higher travel time savings for high-bidding vehicles. 

Several efforts have focused on more microscopic vehicle trajectory optimization formulations to control AVs at traffic intersections. \citet{gregoire2013back} developed a cooperative motion-planning algorithm using a path velocity decomposition to optimally coordinate vehicles while preventing collisions. \citet{de2015autonomous} integrated a priority-based assignment into autonomous intersection simulation. \citet{altche2016analysis} developed a Mixed-Integer Linear Programming (MILP) formulation to coordinate vehicles through intersections. In this model, the IM decides when AVs enter the intersection and at which speed. The authors discretized time and tracked vehicle movement in continuous space variables. \citet{zhang2016optimal,zhang2017decentralized} focused on the decentralized control of traffic intersections based on First-In-First-Out (FIFO) conditions and considered fuel consumption, throughput and turning movements. This framework was extended by \citet{malikopoulos2018decentralized} proposed a decentralized energy-optimal control framework for connected and automated vehicles. Vehicle separation is ensured by rear-end and lateral collision avoidance constraints and the authors prove the existence of a nonempty solution space. 

\citet{levin2017conflict} proposed a conflict point formulation for optimizing intersection throughput referred to as the AIM* protocol. As in \citet{altche2016analysis}, the IM decides on AVs' entry time and speed but instead of discretizing time, the authors discretize the intersection and focus on conflict points to ensure safety. The authors prove that there always exist a conflict-free, feasible solution to the proposed MILP. The present paper builds on this research to coordinate traffic operations during blue phases. 

Only a few papers have addressed the configuration wherein AVs and legacy vehicles share traffic intersections. \citet{dresner2006human} proposed to periodically activate a traditional green phase to allow legacy vehicles to access the intersection. \citet{conde2013intelligent} suggested that legacy vehicles could reserve additional space to ensure safety while using the AIM protocol. \citet{qian2014priority} discussed the use of car-following models for collision avoidance among legacy vehicles and AVs. These efforts rely on the deployment of the AIM protocol to handle vehicle reservations but do not discuss their impact at a network level, \emph{i.e.} the effect of such policies on network throughput and stability. \citet{levin2016multiclass} found that reserving extra space for legacy vehicles created significant delays and that high AV market penetrations (around 80\%) were needed for AIM to improve over traffic signals. The hybrid network traffic control policy combining blue and green phases proposed in this paper may provide an hybrid approach to retain efficiency at lower market penetrations.

\subsection{Our Contributions}

In this paper, we propose a new decentralized, stochastic model for coordinating traffic composed of LVs and AVs in a network of intersections. Our model assumes that vehicle route choice, or equivalently, turning proportions are known. We assume that lane queues are measured periodically and we propose a decentralized, pressure-based algorithm to optimize network throughput. We distinguish between traditional green phases and blue phases which are only available for AVs. For this, we assume that intersections are accessible through two types of lanes: LV-lanes, which can be used by both LVs and AVs, and AV-lanes, which are restricted to AVs. During blue phases, all incoming lanes have access to the intersection and traffic is coordinated using the conflict-point formulation proposed by \citet{levin2017conflict}. At each intersection and each time period, a phase is activated based on the current state of the network using the proposed hybrid pressure-based policy. 

We make the following contributions. 1) We propose a new model for network traffic control with LVs and AVs which combines traditional green phases with AV-restricted blue phases. 2) We present a new MILP formulation for green phase activation wherein turning movement capacity is determined endogenously and lane FIFO blocking effects are accounted for. 3) We extend the max-pressure formulation proposed by \citet{varaiya2013max} to lane-based queues and we propose a new hybrid max-pressure network control policy wherein LVs and AVs share the infrastructure. 4) We characterize the stability region of this system and prove that the proposed hybrid max-pressure network control policy is stable, \emph{i.e.} stabilizes the largest possible set of demand rates in the network. 5) We conduct numerical experiments to test the proposed hybrid network control policy.

The remainder of the paper is organized as follows. We present our stochastic network traffic control model in Section \ref{net}. Intersection phase activation models are presented in Section \ref{phases}. The proposed control policy and its stability proof are introduced in Section \ref{policy}. Numerical results are presented in Section \ref{num} and we discuss our findings in Section \ref{con}.

%%%%%%%%%%%%%%%%%%%%%%%%%%%%%%%%%%%%%%%%%%%%%%%%%%%%%%%%%%%%%%%%%%%%%%%%%%%%
\section{Stochastic Network Traffic Control}\label{net}
%%%%%%%%%%%%%%%%%%%%%%%%%%%%%%%%%%%%%%%%%%%%%%%%%%%%%%%%%%%%%%%%%%%%%%%%%%%%

Consider a traffic network $\G = (\N,\A)$ with a set of intersections $\N$ connected by a set of links $\A$. The set of links is partitioned into three subsets: internal links that connect two intersections, denoted $\Ao \subset \A$, source links at which vehicles enter the network, denoted $\Ar \subset \A$, and sink links at which vehicles exit, denoted $\As \subset \A$. We model congestion at network intersections using point-queues of infinite size and we are interested in the evolution of queue lengths over the entire network. 

We consider two classes of vehicles and lanes: autonomous vehicles (AVs), denoted $a$, and legacy (or human-driven) vehicles (LVs), denoted $l$. We assume that each link of the network consists of a set of lanes which are either restricted to AVs (AV-lanes) or available to both AVs and LVs (LV-lanes). We use $\A_a$ and $\A_l$ to denote AV-lanes and LV-lanes, respectively, and we assume that vehicle movement between different classes of lanes are forbidden. Although AVs can use LV-lanes, we do not model any type of interaction at the traffic-flow level among AVs and LVs on LV-lanes: if an AV uses an LV-lane we assume that it behaves as an LV. We assume that each class of lanes is served by a color-coded traffic phase. Specifically, we assume that LV-lanes are served exclusively by traditional green signal phases. In turn, we assume that AV-lanes are served exclusively by signal-free blue phases. Blue phases differ from green phases in that they directly control vehicles' trajectory within the intersection. The proposed hybrid network control policy presented hereafter chooses which phase (green or blue) should be activated at each intersection of the network and each time period over a discretized time horizon. Both traffic control phases are formally introduced in Section \ref{phases}. \newline

Let $\xit \in \Re_+$ be the number of vehicles on link $i \in \A$ seeking to enter the intersection at time $t$. Although we discretize vehicles in our numerical experiments, integer queue lengths are not necessary for the analytical results presented hereafter. Let $\xt$ be the array of all queue lengths at time $t$. $\xt$ is the state of the network, and the state space is $\X = \{\xit \geq 0 : i \in \A\}$. We consider discretized time and we assume fixed phase time of length $\dt$. Further, we assume that the state of the network $\xt$ is known at each time period $t = 0, \dt, 2\dt, \ldots$. The goal is to design a throughput-optimal network traffic control policy that optimally selects a traffic signal control at each time period $[t, t+\dt[$.

The proposed stochastic network traffic control formulation is based on the concept of vehicle \emph{movements} which are formally defined below.

\begin{defi}
A movement $(i,j) \in \A^2$ is a vehicle trajectory from lane $i$ to lane $j$ across a common intersection $n \in \N$ in the network. We denote $\M$ the set of all movements in the network. 
\end{defi}

AVs communicate their position with IMs to make use of the AIM protocol, hence we assume that movement-specific queues are known for AVs. In contrast, LV-queues can be detected through loop detectors and flow sensors currently in use for traffic signals but their destination is assumed unknown. Specifically, let $\A_a \subset \A$ be the set of AV-restricted lanes. For these lanes, we assume known movements queues, \emph{i.e.} if $j \in \A_a$, then $\xit = \sum_{j \in \A_a} \xijt$ and $\xijt$ is known. In contrast, for other lanes $i \in \A_l = \A \setminus \A_a$, only $\xit$ is known since route choice for LVs is assumed unknown.\\

For each lane $i \in \A$ we assume that lane capacity $C_i$ is known and determined assuming a triangular fundamental diagram relationship, that is $C_i =\frac{\maxU_i w K}{\maxU+w}$, where $\maxU_i$ is the free-flow speed on lane $i$, $K$ is the jam density and $w$ is the congestion wave speed. We also assume that the maximum, unconditional movement service rate $\usij$ is known for each movement and determined based on lost time, specifically: unconditional movement service rates for all movements $(i,j)$ are calculated as $\usij = \min\{C_i, C_j\} \frac{\Delta t - L}{\Delta t}$.

Let $\pijt$ be a random variable denoting the turning proportion from lane $i$ to $j$ at $t$ with known mean $\pij$. Let $\dit$ be a random variable denoting the external incoming traffic onto lane $i$ at time $t$ with known mean $\di$. We denote $\bm{p}$, $\bm{d}$ and $\bm{\bar{s}}$ the vectors of mean turning proportions, mean demands and unconditional movement service rates, respectively. These upper bounds on movement service rates represent the maximum number of vehicles that can be moved from lane $i$ to lane $j$ during time period $t$ when no conflicting movement with $(i,j)$ is activated. From these rates, we can determine maximum, unconditional lane service rates $\usi = \sum_{j \in \A : (i,j) \in \M} \usij$. For convenience if $i$ and $j$ do not correspond to a possible movement in the network, \emph{e.g.} they belong to different intersections or they are both entry or exit lanes of the same intersection, we assume that $\pij = 0$. Hence, we can define the sets of AV-movements $\M_a \equiv \{(i,j) \in \A_a^2 : \pij \neq 0\}$ and LV-movements $\M_l \equiv \{(i,j) \in \A_l^2 : \pij \neq 0\}$. \\

Traffic at network intersections is coordinated by phases which are determined by the selection of the \emph{activation matrix}. In addition, we also introduce the concept of \emph{service matrix} which is used in the proposed traffic control formulations.

\begin{defi}
An activation matrix $\bt$ is a $|\A|\times|\A|$-matrix wherein all entries take value 0 (inactive) or 1 (active), \emph{i.e.} $\bijt \in \{0,1\}$ at time $t$.
\end{defi}
\vspace{-0.5cm}
\begin{defi}
A service matrix $\at$ is a $|\A|\times|\A|$-matrix wherein all entries take a value between 0 (not serviced) and 1 (fully serviced), \emph{i.e.} $\aijt \in [0,1]$ at time $t$.
\end{defi}

The entries of an activation matrix characterize the \emph{activeness} of the corresponding phase: $\bijt = 1$ means that movement $(i,j)$ is active during phase $t$ whereas $\bijt = 0$ means that movement $(i,j)$ is inactive. The entries of the service matrix characterize the \emph{service level} of active movements during phase $t$: $\aijt = 1$ corresponds to a maximal service level for of movement $(i,j)$, whereas $\aijt = 0$ means that movement $(i,j)$ cannot serve any vehicles during time period $t$. Fractional service level values model situations where conflicting movements, \emph{i.e.} posing a safety risk, simultaneously have non-zero activation values. The activation and service matrices are linked through the movement-based constraints $\aijt \leq \bijt$. These linking constraints ensure that an intersection can only service vehicles on movement $(i,j)$ if this movement is active $\bijt = 1$. Further, for priority movements, we set $\aijt = \bijt$ which implies that an activated priority movements has full service level. In the proposed traffic control policy, the selection of the activation matrix requires the solution of two mathematical optimization problems presented in Section \ref{green} and \ref{blue}, respectively and further details are discussed therein.\\

Let $\Yit$ be a random variable with mean $\yit$ denoting the number of vehicles serviced in lane $i$ at time $t$. The vector $\Yt$ is endogenous to the service matrix $\at$ selected by the control policy. Note that the mean $\yit$ of $\Yit$ is unknown and $\yit$ will be modeled as a control variable in the proposed traffic phase optimization formulation. Specifically, the dependency between the vector of lane service rates $\yt$ and $\at$ will be presented in detail in Sections \ref{green} and \ref{blue}, for green and blue phases, respectively. The proposed stochastic network traffic control model is summarized by the lane-queue evolution equation \eqref{eq:queue}.

\begin{equation}
\xjtt = \xjt - \Yjt + \sum_{i \in \A} \pijt \Yit + D_j(t) \quad \forall j \in \A
\label{eq:queue}
\end{equation}

Note, if $j \notin \A_r$ then $\djt = 0$. Conversely, if $j \in \A_r$ then $j$ has no predecessor links, thus $\sum_{i \in \A} \pijt = 0$. Although AV-lanes and other lanes have identical queue evolution equations, the information available is more accurate for AV-lanes, this used in the calculation of network control policies.\\

In related works, the lane service rate vector $\yt$ is commonly calculated as the minimum between the supply (lane or movement service capacity) and the demand (lane or movement queue length) \citep{varaiya2013max,le2015decentralized}. However, this modeling approach does not always capture the interdependency between the activation of possibly conflicting movements with lane or movement service capacity and queue length. Precisely, previous efforts have assumed that a set of activation matrices is provided for each intersection and that lane or movement service capacities can be pre-processed accordingly. This overlooks the impact of queue length on intersection capacity. In contrast our proposed integrated approach aims to accurately estimate the expected number of serviced vehicles at each time period by leveraging the available lane-queue length information (note that this information is also assumed available in the aforementioned papers).

%%%%%%%%%%%%%%%%%%%%%%%%%%%%%%%%%%%%%%%%%%%%%%%%%%%%%%%%%%%%%%%%%%%%%%%%%%%%
\section{Traffic Control with Green and Blue Phases}\label{phases}
%%%%%%%%%%%%%%%%%%%%%%%%%%%%%%%%%%%%%%%%%%%%%%%%%%%%%%%%%%%%%%%%%%%%%%%%%%%%

In this section, we present two intersection-based formulations to coordinate traffic at intersections during green (Section \ref{green}) and blue (Section \ref{blue}) phases, respectively. The proposed  formulations aim to identify the optimal traffic phase at each intersection of the network and at each time period. The notion of cycles, i.e. a sequence of phases, is not represented in the proposed network model. Specifically, we do not impose any constraints on the sequence of phases that may be activated at network intersections. Instead, this is left to be optimized by the proposed network traffic control policy, which we will prove to be throughput-optimal (Section \ref{policy}). Notably, we show that network-wide traffic coordination can be decentralized by implementing a hybrid pressure-based policy that activates the highest pressure phase---among green and blue phases---at each intersection and at each time period.

%%%%%%%%%%%%%%%%%%%%%%%%%%%%%%%%%%%%%%%%%%%%%%%%%%%%%%%%%%%%%%%%%%%%%%%%%%%%
\subsection{Green Phase}\label{green}
%%%%%%%%%%%%%%%%%%%%%%%%%%%%%%%%%%%%%%%%%%%%%%%%%%%%%%%%%%%%%%%%%%%%%%%%%%%%
For green phases, we propose a MILP approach that aims to identify the optimal activation matrix $\at$ while accounting for movement capacity loss due to potential conflicting movements and known lane-queues. Specifically, we estimate the expected number of serviced vehicles (or lane service rates) for each lane in the network when this value is assumed endogenous to the selected activation matrix. 

Let $\A_l^n \subset \A_l$ (respectively, $\M_l^n \subset \M_l$) be the set of LV-lanes (respectively, LV-movements) of intersection $n \in \N$. We make the following design assumptions:

\begin{enumerate}
\item The set of LV-movements $\M_l^n$ can be partitioned into two sets: priority $\P^n$ and yield  $\Y^n$ movements, \emph{i.e.}  $\M_l^n = \P^n \cup \Y^n$ and $\P^n \cap \Y^n = \emptyset$.
\item Two conflicting priority movements cannot be selected simultaneously.
\item Two conflicting yield movements cannot be selected simultaneously.
\end{enumerate}

Assumptions 2 and 3 are motivated by the fact that in existing traffic intersections there needs to be a consensus between movements posing a safety risk: this consensus is traditionally resolved by ``right-of-way'' rules which implies that one movement must have priority over the other. For instance, two through movements cannot be activated simultaneously unless they are parallel. In our numerical experiments we model through and right turns as priority movements whereas left turns are categorized as yield movements. 

To implement the intersection design assumptions listed above, we introduce two binary matrices that can be pre-processed for all intersections in the network. Formally, let $\C_{ij}$ be the set of conflict points that movement $(i,j) \in \M_l^n$ passes through. This set can be defined by examining the geometry of the intersection containing movement $(i,j)$ and identifying conflict points with other movements---an illustration is provided in Section \ref{example}. Given conflict point sets $\C_{ij}$ for each movement $(i,j) \in \M^n_l$, we can determine the values $\cijij$ for each pair of movements as follows:
\begin{equation}
\cijij = \begin{cases}
1 \text{ if } \C_{ij} \cap \C_{i'j'} \neq \emptyset \\
0 \text{ otherwise}
\end{cases}\quad \forall (i,j), (i',j') \in \M_l^n
\label{eq:conflicts}
\end{equation}

Let $g_{ij}^{i'j'}$ be a binary parameter indicating forbidden simultaneous movements, formally defined as:
\begin{equation}
\gijij = \begin{cases}
1 \text{ if } (i,j), (i',j') \in \P \text{ and } \cijij = 1 \\
1 \text{ if } (i,j), (i',j') \in \Y \text{ and } \cijij = 1 \\
0 \text{ otherwise}
\end{cases}\quad \forall (i,j), (i',j') \in \M_l^n
\label{eq:forbidden}
\end{equation}

Let $\bijt \in \{0,1\}$ be the binary variable representing the activation of LV-movement $(i,j) \in \M_l^n$ in the activation matrix ($\bijt = 1$) or not ($\bijt = 0$). To forbid the simultaneous activation of two priority or two yield conflicting movements, as required by the design assumptions 2 and 3, we impose the following constraint: 
\begin{equation}
\bijt + \bijpt \leq 1 \quad \forall (i,j), (i',j') \in \M_l^n : g_{ij}^{i'j'} = 1
\label{eq:forbid}
\end{equation}

We can use a supply-demand formulation to calculate the lane service rates $\yijt$ for movements $(i,j) \in \M_l^n$. Let $\aijt \in [0,1]$ be the decision variable representing the fraction of capacity allocation to LV-movement $(i,j) \in \M_l^n$. Recall that We have the relationship $\aijt \leq \bijt$. On the supply side, the expected movement capacity is $\aijt \usij$. Assuming known, exogenous turning proportions, $\pij$, the demand for movement $(i,j)$ at time $t$ is upper-bounded by $\pij \xit$. However, FIFO behavior on lanes may lead to vehicle-blocking \citep{tampere2011generic,li2016effects}. To capture this lane FIFO behavior, let $\pit \in [0,1]$ be a variable representing FIFO blocking effects on vehicles movements from lane $i \in \A_l^n$, i.e. $\pit = 1$ means no blocking effects, whereas $\pit < 1$ implies that only a fraction of the demand for each turning movement is able to traverse the intersection due to FIFO blocking effects. For instance, if only 60\% of the demand for a given movement can exit its upstream lane, then the remaining 40\% blocks the demand of other turning movements. Hence, in this example, only 60\% of the demand for other turning movements may exit and $\pit = 0.6$. The expected service rate of movement $(i,j) \in \M_l^n$ can then be defined as:
\begin{equation}
\yijt \equiv \min\{\aijt \usij, \pij \xit \pit\}
\label{eq:service}
\end{equation}

The value of lane-based variables $\pit$ can be adjusted based on the \textit{supply-to-demand ratio} of each movement. For each lane $i$, we define the level of FIFO blocking effects as the minimum supply-to-demand ratio and the no-blocking case, i.e. $\pit=1$. This yields the following definition:
\begin{equation}
\pit \equiv \min\left\{\min_{(i,j) \in \M_l^n}\left\{\frac{\aijt \usij}{\pij \xit} \right\}, 1 \right\} \quad \forall i \in \A_l^n
\label{eq:fifo}
\end{equation}

Observe that if $\min_{(i,j) \in \M_l^n}\left\{\frac{\aijt \usij}{\pij \xit} \right\} \geq 1$ for all movements $(i,j) \in \M_l^n$ emanating of lane $i \in \A_l^n$, then $\pit = 1$ and $\aijt \usij \geq \pij \xit = \pij \xit \pit$, thus $\yijt = \pij \xit \pit$ for all movements $(i,j) \in \M_l^n$. Otherwise, if there exists an outgoing lane $j$ such that $\frac{\aijt \usij}{\pij \xit} < 1$ and this ratio is minimal over all emanating movements from lane $i$; then $\pit = \frac{\aijt \usij}{\pij \xit}$ and $\yijt = \aijt \usij = \pij \xit \pit$. Hence, the FIFO condition \eqref{eq:fifo} implies that the expected service rate $\yijt$ of movement $(i,j) \in \M_l^n$ defined in \eqref{eq:service} is always demand-constrained, i.e. $\yijt = \pij \xit \pit$. Movement service rates can be linked to lane service rates by summing over all movements from the incoming lane. Since movement service rates are always demand-constrained, we obtain $\yit = \sum_{j \in \A_l^n : (i,j) \in \M_l^n} \yijt = \sum_{j \in \A_l^n : (i,j) \in \M_l^n}\pij\xit\pit = \xit\pit$, for all LV lanes $i \in \A_l^n$.\newline

For a priority movement $(i,j) \in \P^n$, recall that $\aijt = \bijt$. In turn, for a yield movement $(i,j) \in \Y^n$, the available capacity depends on the selection of conflicting movements. Note that, based on our design assumptions, a yield movement may conflict with any number of priority movements as long as these priority movements do not conflict with each other. To calculate the endogenous expected capacity of yield movements, we define $\mij \geq 0$ as the \emph{slack} of movement $(i,j) \in \M_l^n$ representing the expected available capacity for this movement if the activation matrix $\bt$ is selected. Hence, the slack of movement $(i,j) \in \M_l^n$ is equivalent to the reserve capacity of movement $(i,j)$ if $(i,j)$ is activated. Mathematically, the slack of $(i,j) \in \M_l^n$ can be defined as:
\begin{equation}
\mij \equiv \left(\usij - \yijt\right)\bijt = \left(\usij - \pij \xit \pit\right)\bijt
\label{eq:slack}
\end{equation}

Recall that $\C_{ij}$ is the set of conflicts movements of $(i,j) \in \M_l^n$. The endogenous expected capacity of yield movements depends on the choice of the activation matrix. The available capacity of a yield movement $(i,j) \in \Y^n$ is defined as the minimum slack of its conflicting activated movements. This definition ensures that the available capacity of $(i,j)$ is no more than the minimum slack available among active movements in $\C_{ij}$, while inactive movements in $\C_{ij}$ do not restrict the available capacity of $(i,j)$. Hence, we can define the expected capacity of yield movement $(i,j) \in \Y^n$ as: 
\begin{equation}
\aijt \usij = \min\left\{\min_{(i'j') \in \C_{ij}}\left\{\mijp : \bijpt = 1\right\}, \usij \bijt\right\}
\label{eq:endo}
\end{equation}

To incorporate Equations \eqref{eq:fifo}, \eqref{eq:slack} and \eqref{eq:endo} as MILP constraints, we linearize these equations exactly using additional binary variables. Let $\chit \in \{0,1\}$ and $\chijt \in \{0,1\}$ be binary variables used to linearize the right-hand side of \eqref{eq:fifo}. Equation \eqref{eq:fifo} is equivalent to the constraints:
\begin{subequations}
\begin{align}	
& && \pit \leq 1 && \forall i \in \A_l^n \\ 	
& && \pit \geq \chit && \forall i \in \A_l^n \\ 	
& && \pit \leq \frac{\aijt \usij}{\pij \xit} && \forall (i,j) \in \M_l^n : \pij \xit > 0 \\	
& && \pit \geq \frac{\aijt \usij}{\pij \xit} - (1 -  \chijt)\frac{\usij}{\pij \xit}&& \forall (i,j) \in \M_l^n : \pij \xit > 0 \\
& && \chit + \sum_{j : (i,j) \in \M_l^n, \pij \xit > 0} \chijt = 1 && \forall i \in \A_l^n \\
& && \chit \in \{0,1\} && \forall i \in \A_l^n \\ 		
& && \chijt \in \{0,1\} && \forall (i,j) \in \M_l^n 	
\end{align}
\label{eq:linpit}
\end{subequations}

To exactly linearize Constraints \eqref{eq:slack} and \eqref{eq:endo}, let $\lijt \in \{0,1\}$ and $\lijpt \in \{0,1\}$ be binary variables used to model the right-hand side of \eqref{eq:endo}. Equations \eqref{eq:slack} and \eqref{eq:endo} are equivalent to the constraints:
\begin{subequations}
\begin{align}
& && \aijt \leq \bijt && \forall (i,j) \in \Y^n \\
& && \aijt \geq \bijt - (1 - \lijt) && \forall (i,j) \in \Y^n \\	
& && \aijt \usij \leq \bijpt \usijp - \pijp \xipt \pipt \nonumber \\
& &&  + (1 - \bijpt)\usij && \forall (i,j) \in \Y^n, (i',j') \in \C_{ij} \\
& && \aijt \usij \geq \bijpt \usijp - \pijp \xipt \pipt \nonumber \\
& && - (1 - \bijpt)\usij - (1  - \lijpt)\max\{\usij,\usijp\} && \forall (i,j) \in \Y^n, (i',j') \in \C_{ij} \\	
& && \lijt + \sum_{(i',j') \in \C_{ij}} \lijpt = 1 && \forall (i,j) \in \Y^n \\
& && \lijpt \leq \bijpt && \forall (i,j) \in \Y^n, (i',j') \in \C_{ij} \\	
& && \lijt \in \{0,1\} && \forall (i,j) \in \Y^n \\ 	
& && \lijpt \in \{0,1\} && \forall (i,j) \in \Y^n, (i',j') \in \C_{ij} 	
\end{align}
\label{eq:linyield}
\end{subequations}

Our objective is to select a control policy that maximizes network throughput. We build on the max-pressure literature and define $\wit$ as the \emph{weight} of lane $i \in \A$ at time $t$ onto the network based on current queue lengths \citep{tassiulas1992stability}:
\begin{equation}
\wit = \xit - \sum_{j \in \A : (i,j) \in \M} \pij \xjt   
\label{eq:weight}
\end{equation} 

As in \citet{varaiya2013max}, we will later show that maximizing pressure locally, \emph{i.e.} at each intersection, will maximize network throughput. Hence, the objective function for the green phase activation program should maximize local pressure, \emph{i.e.} $\sum_{i \in \A_l^n} \wit \yit = \sum_{i \in \A_l^n} \wit \xit \pit$.

Let $\zgreen$ be the maximal local pressure that can be obtained using the green phase based on the network state $\xt$ at intersection $n \in \N$. The proposed mathematical programming formulation for identifying maximal pressure green phases is summarized below in \eqref{mod:green} and hereby referred to as the \green.

\begin{subequations}
\begin{align}
	\zgreen =\ & \max && \sum_{i \in \A_l^n} \wit \xit \pit \label{eq:objgreen}\\
& \mathrm{s.t.}  && \pit \leq 1 && \forall i \in \A_l^n \\ 	
& && \pit \geq \chit && \forall i \in \A_l^n \\ 	
& && \pit \leq \frac{\aijt \usij}{\pij \xit} && \forall (i,j) \in \M_l^n : \pij \xit > 0 \\	
& && \pit \geq \frac{\aijt \usij}{\pij \xit} - (1 -  \chijt)\frac{\usij}{\pij \xit}&& \forall (i,j) \in \M_l^n : \pij \xit > 0 \\
& && \chit + \sum_{j : (i,j) \in \M_l^n, \pij \xit > 0} \chijt = 1 && \forall i \in \A_l^n \\ 			
& && \bijt + \bijpt \leq 1 && \forall (i,j), (i',j') \in \M_l^n : g_{ij}^{i'j'} = 1 \\
& && \aijt = \bijt && \forall (i,j) \in \P^n \\
& && \aijt \leq \bijt && \forall (i,j) \in \Y^n \\
& && \aijt \geq \bijt - (1 - \lijt) && \forall (i,j) \in \Y^n \\	
& && \aijt \usij \leq \bijpt \usijp - \pijp \xipt \pipt \nonumber \\
& &&  + (1 - \bijpt)\usij && \forall (i,j) \in \Y^n, (i',j') \in \C_{ij} \\
& && \aijt \usij \geq \bijpt \usijp - \pijp \xipt \pipt \nonumber \\
& && - (1 - \bijpt)\usij - (1  - \lijpt)\max\{\usij,\usijp\} && \forall (i,j) \in \Y^n, (i',j') \in \C_{ij} \\	
& && \lijt + \sum_{(i',j') \in \C_{ij}} \lijpt = 1 && \forall (i,j) \in \Y^n \\
& && \lijpt \leq \bijpt && \forall (i,j) \in \Y^n, (i',j') \in \C_{ij} \\	
& && \chit \in \{0,1\} && \forall i \in \A_l^n \\ 		
& && \chijt \in \{0,1\} && \forall (i,j) \in \M_l^n \\ 			
& && \lijt \in \{0,1\} && \forall (i,j) \in \Y^n \\ 	
& && \lijpt \in \{0,1\} && \forall (i,j) \in \Y^n, (i',j') \in \C_{ij} \\ 		
& && \bijt \in \{0,1\} && \forall (i,j) \in \M_l^n \\ 
& && \aijt \geq 0 && \forall (i,j) \in \M_l^n \\ 
& && \pit \geq 0 && \forall i \in \A_l^n 
\end{align}
\label{mod:green}
\end{subequations}

\begin{figure}[ht]
	\centering
	\includegraphics[width=0.7\linewidth]{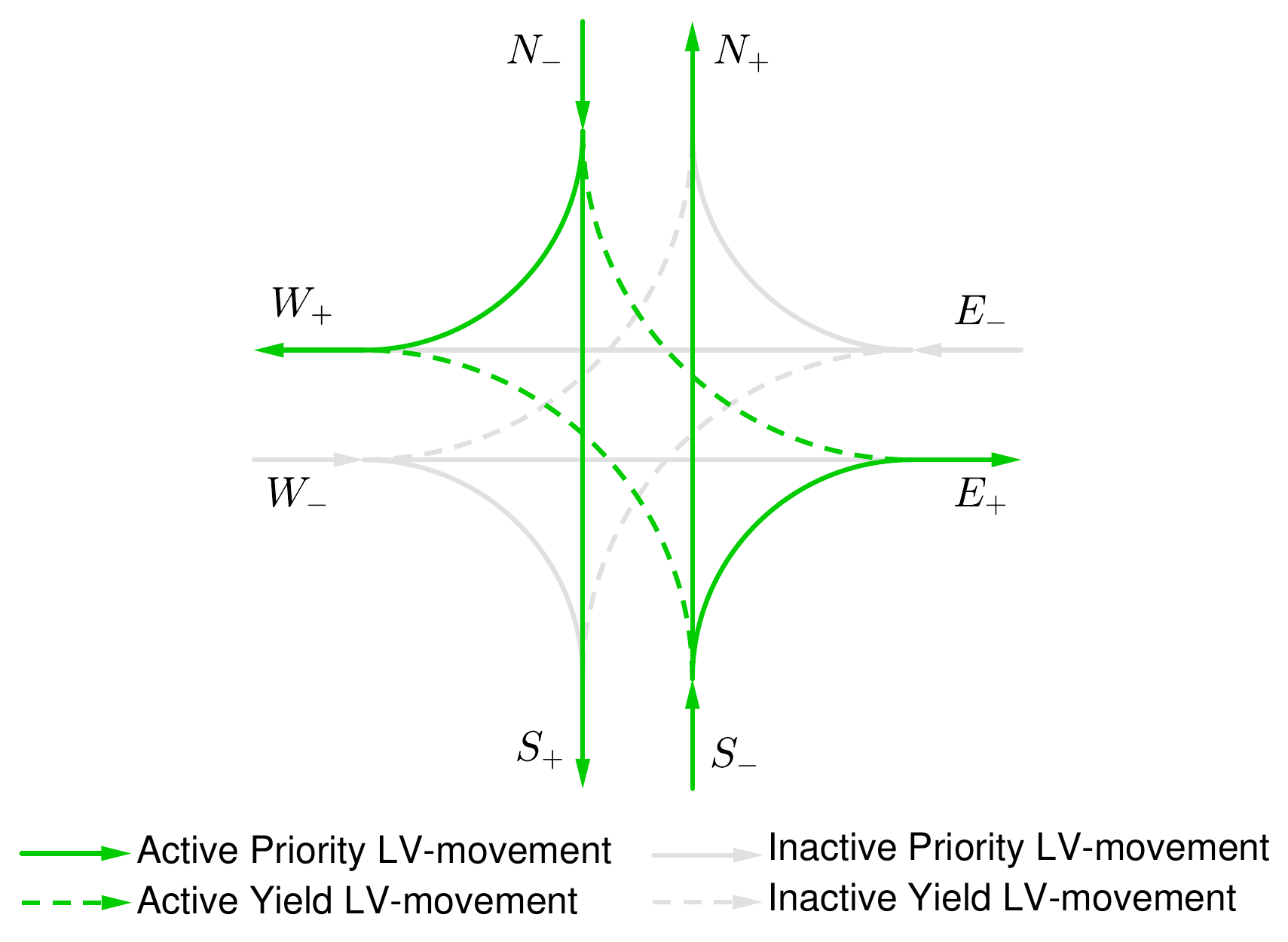}
	\caption{Illustration of a green phase activation.}
	\label{fig:greenphase}
\end{figure}

To illustrate the behavior of \green, consider a 4-approach intersection where each lane has 3 movements (left, straight, right). All straight and right turns are priority movements whereas left turns are yield movements. Assume that all movements from South and North lanes are activated as illustrated in Figure \ref{fig:greenphase}. We examine the behavior of the model from the perspective of lane $S_-$. Without any loss of generality, we assume that all movements have identical capacity $\usij$, denoted $\bar{s}$ for short. For this activation matrix, the conflict set of yield movement $(S_-,W_+)$ is the opposite priority movements: $\C_{S_-W_+} =\{(N_-,S_+),(N_-,W_+)\}$. Observe that Constraints \eqref{eq:linyield} for movement $(S_-,W_+)$ collapse to:
\begin{align*} 
\alpha_{S_-W_+} & \leq 1 \\
\alpha_{S_-W_+} & \geq 1 - (1 - \lambda_{S_-W_+}) \\
\alpha_{S_-W_+} \bar{s} & \leq \bar{s} - p_{N_-S_+} x_{N_-}\phi_{N_-} \\
\alpha_{S_-W_+} \bar{s} & \leq \bar{s} - p_{N_-W_+} x_{N_-}\phi_{N_-} \\
\alpha_{S_-W_+} \bar{s} & \geq \bar{s} - p_{N_-S_+} x_{N_-}\phi_{N_-} - (1 - \lambda_{S_-W_+}^{N_-S_+})\bar{s} \\
\alpha_{S_-W_+} \bar{s} & \geq \bar{s} - p_{N_-W_+} x_{N_-}\phi_{N_-} - (1 - \lambda_{S_-W_+}^{N_-W_+})\bar{s}\\
\lambda_{S_-W_+} + \lambda_{S_-W_+}^{N_-S_+} + \lambda_{S_-W_+}^{N_-W_+} &= 1
\end{align*}

If the slack of both movements $(N_-,S_+)$ and $(N_-,W_+)$ is greater than the capacity of $(S_-,W_+)$ (i.e. $\bar{s}$), then $\lambda_{S_-W_+}=1$ and $\alpha_{S_-W_+}=1$. In this case, yield movement $(S_-,W_+)$ is unrestricted. Otherwise, assume that the slack of $(N_-,S_+)$ is lower than that of $(N_-,W_+)$ and lower than $\bar{s}$, then $\lambda_{S_-W_+}^{N_-S_+}=1$ and $\alpha_{S_-W_+}=1 - \frac{p_{N_-S_+} x_{N_-}\phi_{N_-}}{\bar{s}}$. Observe that $\alpha_{S_-W_+} \geq 0$ since Constraints \eqref{eq:linpit} require $p_{N_-S_+} x_{N_-}\phi_{N_-} \leq \alpha_{N_-S_+} \bar{s} = \bar{s}$ since $(N_-,S_+)$ is an activated priority movement. In this case, movement $(N_-,S_+)$ restricts the service level of yield movement $(S_-,W_+)$.

The level of FIFO blocking effects on lane $S_-$ (i.e. $\phi_{S_-}$) takes the minimum value among supply-to-demand ratio of all movements emanating from this lane and 1. If all supply-to-demand ratios are greater than 1, then $\chi_{S_-}=1$ and $\phi_{S_-}=1$, which corresponds to the case where no FIFO blocking effects occur on lane $S_-$. Otherwise, in the continuity of this illustration, assume that the lowest supply-to-demand ratio is achieved for movement $(S_-,W_+)$ (due for instance to the restricted service level of this movement) and this supply-to-demand ratio is lesser than 1. In this case $\chi_{S_-W_+} = 1$ and $\phi_{S_-} = \frac{\alpha_{S_-W_+}\bar{s}}{p_{S_-W_+} x_{S_-}}$. In this case, FIFO blocking effects due to the limited available capacity of movement $(S_-,W_+)$ are restricting the outflow of lane $S_-$, which in turn affect the pressure release of this lane. This phenomena is observed in current lanes where a left-turning vehicle blocks through movements for the lane. In practice, this is often alleviated by adding a left-turn bay.

This analysis highlights that Model \green is fully determined by the choice of the activation matrix $\bt$ which is decided based on the maximum pressure release that can be achieved. The solution of \green gives an optimal service matrix at intersection $n$ for the green phase, denoted $\atng$ and an optimal (binary) activation matrix denoted $\btng$. 

%%%%%%%%%%%%%%%%%%%%%%%%%%%%%%%%%%%%%%%%%%%%%%%%%%%%%%%%%%%%%%%%%%%%%%%%%%%%
\subsection{Blue Phase}\label{blue}
%%%%%%%%%%%%%%%%%%%%%%%%%%%%%%%%%%%%%%%%%%%%%%%%%%%%%%%%%%%%%%%%%%%%%%%%%%%%

To coordinate traffic during blue phases, we adapt a mixed-integer programming formulation from \citet{levin2017conflict} to maximize local pressure. For blue phases, lane service rates can be obtained by solving an optimization problem wherein collision avoidance constraints are imposed at all conflicts point of the intersection. In contrast to \green, the blue phase model optimizes individual vehicle trajectories while ensuring traffic safety. Specifically, the blue phase Model finds optimal AVs speeds and departure time based on the current AV demand at each intersection.\\

Let $\Vt$ be set of AVs in the network at time $t$ seeking to enter intersection $n \in \N$. Let $\A_a^n \subset \A_a$ (respectively, $\M_a^n \subset \M_a$) be the set of AV-lanes (respectively, AV-movements) of intersection $n \in \N$. For each intersection, the set of possible AV-movements is assumed known and intersecting AV-movements generate \emph{conflict-points}. Thus, an AV-movement $(i,j) \in \M_a^n$ can be viewed as a two-dimensional trajectory which consists of a sequence of conflict-points starting at the head node of lane $i$ denoted $i^+$ and ending at the tail node of lane $j$ denoted $j^-$. Since AVs' route choice is assumed to be known, we can map each vehicle to a trajectory. Let $\inc_v$ (respectively, $\out_v$) be the entry (respectively, exit) point of vehicle $v \in \Vt$ into the intersection. Let $\rho_v$ be the trajectory of $v$, \emph{i.e.} $\rho_v = \{\inc_v, \ldots ,\out_v\}$.\\

Recall that the current time period consists of the interval $[t, t + \dt[$. Let $t_v(c) \geq t$ be a decision variable representing the arrival time of vehicle $v$ at point $c \in \rho_v$ and let $\tau_v(c) \geq 0$ be a decision variable representing the time conflict point $c \in \rho_v$ is reserved for vehicle $v$. The values of these variables is determined by the speed assigned and the departure assigned to $v$ as described in \citet{levin2017conflict}.\\

Let $z_v$ be a binary variable denoting if vehicle $v \in \V^n(t)$ traverses intersection $n$ during the current time period (1) or not (0). Formally, 
\begin{equation}
z_v = 1 \quad \Leftarrow \quad t_v(\out_v) + \tau_v(\out_v) \leq t + \dt
\label{eq:bin}
\end{equation}

A relaxed form of this relationship can be modeled using integer-linear constraints as follows:
\begin{equation}
t_v(\out_v) + \tau_v(\out_v) \leq t + \dt + (1 - z_v)M_v
\label{eq:zt}
\end{equation}

where $M_v \geq 0$ represents a large number which can be set to the maximal exit time of vehicle $v \in \Vt$. Constraint \eqref{eq:zt} imposes that $z_v = 1$ if $t_v(\out_v) + \tau_v(\out_v) \leq t + \dt$ and is free otherwise. To ensure that $z$-variables are only activated if all predecessor vehicles in the same queue have traversed the intersection, we impose the constraints $z_{v'} \leq z_v$, for all vehicles $v, v' \in \Vt : \inc_v = \inc_{v'}, e_v < e_{v'}$. Let $\Vit = \{v \in \Vt : \inc_v = i^+\}$ be the set of vehicles seeking to travel from lane $i \in \A_a^n$ at time $t$. The number of AVs serviced on lane $i \in \A_a^n$ is $\sum_{v \in \Vit} z_v$. Similarly to the objective function of \green, the objective of the blue phase model is to maximize local pressure on AV-lanes. This can be formulated as: $\sum_{i \in \A_a^n} \wit \sum_{v \in \Vit} z_v$. 

The remaining of the blue phase model is identical to that presented in \citet{levin2017conflict}. Note that movement capacity is only implicitly represented within the blue phase model. Recall that we assume a triangular fundamental diagram relationship in the intersection. This relationship is used to determine the reservation time $\tau_v(c)$ in Constraint \eqref{eq:tau}. Constraints \eqref{eq:speedbounds} and \eqref{eq:ctespeed} are used to impose lower ($\minU_v$) and upper ($\maxU_v$) bounds on the speed of vehicle $v$, and to impose that vehicles have a constant speed throughout the intersection. FIFO conditions at all shared conflict points for vehicles in the same lane queue are imposed by Constraints \eqref{eq:fifoblue}. Binary variables $\dvvp(c)$ are used to model the order of vehicles at conflict points, \emph{i.e.} $\dvvp(c) = 1$ (respectively $\dvpv(c) = 1$) means that vehicle $v$ (respectively, $v'$) traverses conflict point $c$ before vehicle $v'$ (respectively, $v$) and disjunctive trajectory separation constraints (see \eqref{eq:dis1} and \eqref{eq:dis2} in the formulation below) are used to ensure that conflict points are reserved a sufficient amount of time to ensure traffic safety. This formulation builds on space-discretized collision avoidance formulations for air traffic control \citep{rey2015equity,rey2015subliminal}. For more details on this conflict-point formulation, we refer the reader to \citet{levin2017conflict}. Let $\zblue$ be the maximal local pressure that can be obtained using the blue phase based on the network state $\xt$ at intersection $n$. The MILP used to coordinate traffic during blue phases is summarized below in Formulation \eqref{mod:blue} and hereby referred to as the \blue.

\begin{subequations}
\begin{align}
\zblue =\ & \max && \sum_{i \in \A_a^n} \wit \sum_{v \in \Vit} z_v  \\
& \mathrm{s.t.} && t_v(\out_v) + \tau_v(\out_v) \leq t + \dt + (1 - z_v)M_v && \forall v \in \Vt \\
%& && z_{v'} \leq z_v \quad && \forall v, v' \in \V^n(t) : \inc_v = \inc_{v'}, e_v < e_{v'} \\ we have not implemented this constraint
& && t_v(\inc_v)\geq e_v && \forall v \in \Vt \\
%& && t_v(\inc_v) +\tau_v(\inc_v) \leq  t_{v'}(\inc_{v'}) && \forall v, v' \in \V^n(t) : \inc_v=\inc_{v'}, e_v < e_{v'} \\
& && \tau_v(c) = \frac{L_v}{w} + \frac{L_v(t_v(\out_v)-t_v(\inc_v))}{d_v(\inc_v,\out_v)} && \forall v \in \Vt, \forall c\in\p_v \label{eq:tau}\\
& && \frac{d_v(\inc_v,\out_v)}{\maxU_v} \leq t_v(\out_v) - t_v(\inc_v) \leq \frac{d_v(\inc_v,\out_v)}{\minU_v} && \forall v \in \Vt \label{eq:speedbounds} \\
& && \frac{t_v(c)-t_v(\inc_v)}{d_v(\inc_v,c)} = \frac{t_v(\out_v)-t_v(\inc_v)}{d(\inc_v,\out_v)} && \forall v \in \Vt, \forall c\in\p_v \label{eq:ctespeed} \\
& && t_v(c) + \tau_v(c) \leq t_{v'}(c) && \forall v, v' \in \Vt : \inc_v = \inc_{v'}, e_v < e_{v'}, \nonumber \\
& && && \forall c\in\p_v\cap\p_{v'} \label{eq:fifoblue} \\
& && t_v(c)+\tau_v(c)-t_{v'}(c) \leq (1 - \dvvp(c))M_{vv'} && \forall v, v' \in \V^n(t): \inc_v\neq\inc_{v'}, \nonumber \\ 
& && && \forall c\in\p_v\cap\p_{v'} \label{eq:dis1}\\
& && \dvvp(c)+\dvpv(c)= 1 && \forall v, v' \in \V^n(t): \inc_v\neq\inc_{v'}, v<{v'}, \nonumber \\
& && && \forall c\in\p_v\cap\p_{v'} \label{eq:dis2}\\
& && z_{v'} \leq z_v && \forall v, v' \in \Vt : \inc_v = \inc_{v'}, e_v < e_{v'} \\
& && \dvvp(c) \in \{0,1\} && \forall v, v' \in \V^n(t): \inc_v\neq\inc_{v'}, \nonumber \\
& && && \forall c\in\p_v\cap\p_{v'} \\
& && z_v \in \{0,1\} && \forall v \in \V^n(t) \\
& && t_v(c)\geq t && \forall v \in \V^n(t), \forall c \in \p_v \\
& && \tau_v(c)\geq 0 && \forall v \in \V^n(t), \forall c \in \p_v 
\end{align}
\label{mod:blue}
\end{subequations}

To derive the service matrix $\atnb$ associated with the optimal solution of \blue at intersection $n$ and time $t$, we introduce movement-based vehicle sets $\Vijt = \{v \in \Vt : \inc_v = i^+, \out_v = j^-\}$ for each $(i,j) \in \M_a^n$. The service rate of movement $(i,j) \in \M_a^n$ is then calculated as the ratio of the number of vehicles serviced by the unconditional movement service rate, \emph{i.e.}:
\begin{equation}
\aijt = \frac{\sum_{v \in \Vijt}{z_v}}{\usij}
\label{eq:aijblue}
\end{equation}

The (binary) activation matrix $\btnb$ associated with the optimal solution of \blue intersection $n$ and time $t$ can be defined based on $\atnb$ as follows:
\begin{equation}
\bijt = \begin{cases}
1 \text{ if } \aijt > 0 \\
0 \text{ otherwise}
\end{cases}
\label{eq:bijblue}
\end{equation}

Similarly to green phase, we can also derive the FIFO blocking effects $\pit$ associated to lane $i$ in blue phase using the relationship $\yit = \xit \pit$:
\begin{equation}	
\pit = \frac{\yit}{\xit} = \frac{\sum_{v \in \V^n_i(t)} z_v}{\xit} 
\label{eq:phiblue}
\end{equation}

Although variables $\aijt$, $\bijt$ and $\pit$ are only implicitly defined in \blue, they are used to characterize the stability properties of the proposed network control policy, as discussed in Section \ref{policy}.

%%%%%%%%%%%%%%%%%%%%%%%%%%%%%%%%%%%%%%%%%%%%%%%%%%%%%%%%%%%%%%%%%%%%%%%%%%%%
\subsection{Illustration of \green and \blue Traffic Control Formulations}\label{example}
%%%%%%%%%%%%%%%%%%%%%%%%%%%%%%%%%%%%%%%%%%%%%%%%%%%%%%%%%%%%%%%%%%%%%%%%%%%%

\begin{figure}[ht]
\centering
\begin{subfigure}[b]{0.33\linewidth}
\centering\includegraphics[width=0.95\linewidth]{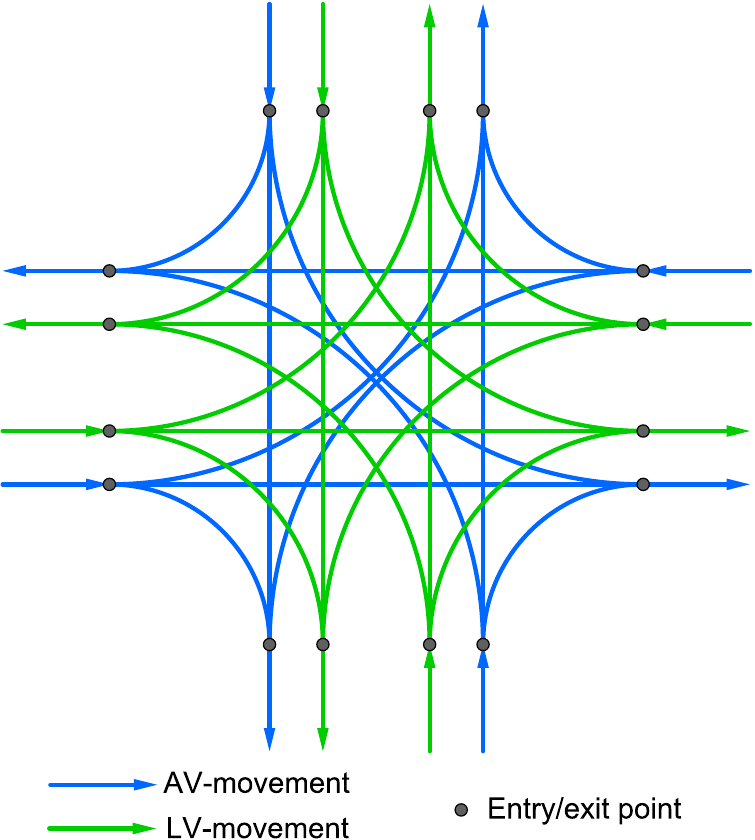}
\caption{\label{fig:interall}}
\end{subfigure}%
\begin{subfigure}[b]{0.33\linewidth}
\centering\includegraphics[width=0.95\linewidth]{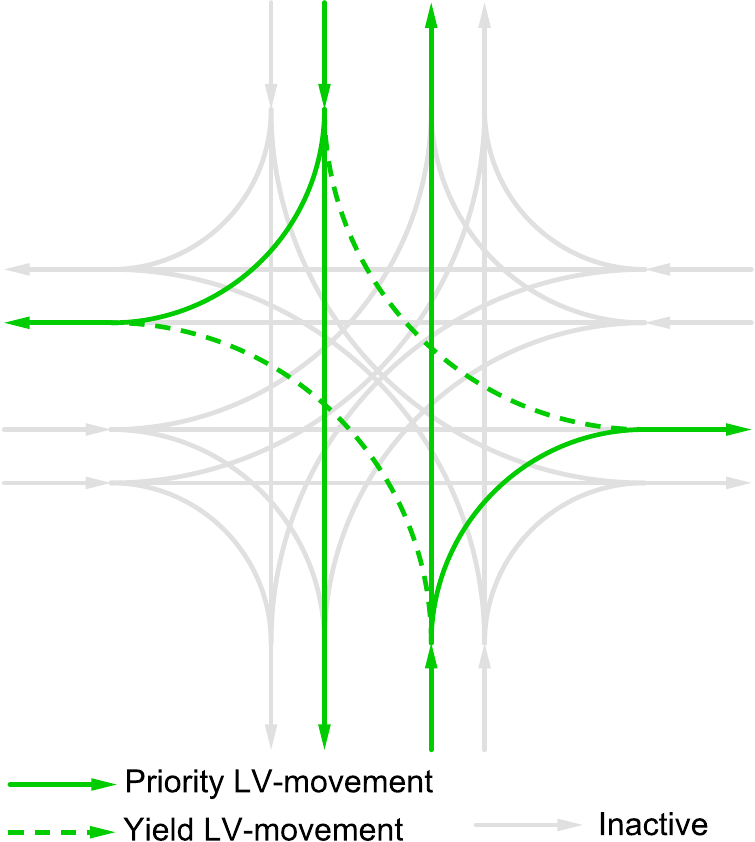}
\caption{\label{fig:interpy}}
\end{subfigure}
\begin{subfigure}[b]{0.33\linewidth}
\centering\includegraphics[width=0.95\linewidth]{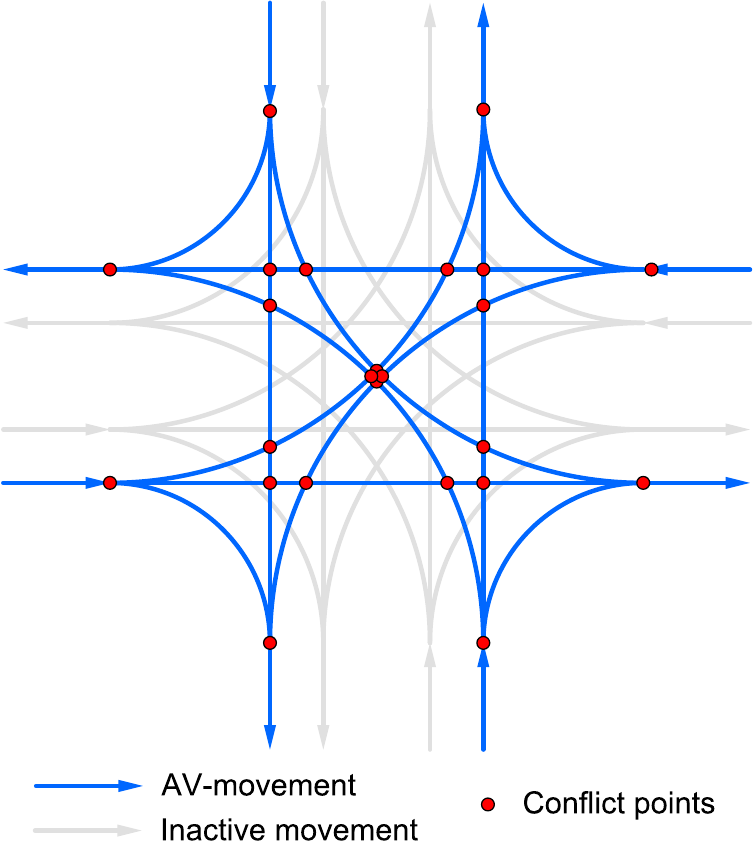}
\caption{\label{fig:interblue}}
\end{subfigure}	
\caption{Figure~\ref{fig:interall} shows all possible LV- and AV-movements on a typical intersection. Figure~\ref{fig:interpy} depicts a possible green phase involving priority and yield LV-movements and Figure~\ref{fig:interblue} illustrates the conflict-point formulation during blue phase.\label{fig:inter}}
\end{figure}

The proposed traffic control formulation is illustrated in Figure \ref{fig:inter} for a typical traffic intersection connected to eight incoming and eight outgoing links, each of which composed of one LV-lane and one AV-lane. A possible solution to \green is shown in Figure \ref{fig:interpy} with four priority LV-movements and two yield LV-movements. The blue phase is illustrated in \ref{fig:interblue} wherein AV-traffic is coordinated by solving \blue and reserving conflicts points for each vehicle in order to maximize local pressure. 

To further illustrate the behavior of formulations \green and \blue, we examine a specific instance based on the intersection geometry illustrated in Figures~\ref{fig:interpy} and~\ref{fig:interblue}, respectively. We consider a 4-approach intersection with one LV-lane and one AV-lane per approach. Lanes are labelled based on their cardinal orientation. Incoming lanes are sub-scripted with a $-$ sign and outgoing lanes are sub-scripted with a $+$ sign. Each lane has three possible movements: right, through and left. For green phase purposes, Right and Through movements are assumed to be \emph{priority} ($\P$) movements whereas left movements are assumed to be \emph{yield} ($\Y$) movements. 

We assume that the free-flow speed is $\maxU = 44$ ft/s and the congestion wave speed is $w = 11$ ft/s. We consider an intersection of width 48 ft and assume all vehicles to be of length 17.6 ft. We use a time period of $\Delta t = 10$ s and we assume a lost time of 2 s for green phase which represents all-red clearance time and vehicle start-up delay. Under this configuration, lane capacity ($C_i$) is 5 vehicles per time period. The unconditional movement service rate of green phase is $\usij^G = 4$ vehicles per time period, and, assuming no lost time for blue phase, the unconditional movement service rate of blue phase is $\usij^B = 5$ vehicles per time period. This information, along with turning proportions and movement specific information is detailed in Table \ref{tab:mov}.

\begin{table}[b]%
\centering
\begin{tabular}{lllllll}
\toprule
& & & $C_i$ & $\usij^G$ & $\usij^B$ & Number of \\ 
Movement & Type & $\pij$ & (in veh/$\Delta t$) & (in veh/$\Delta t$) & (in veh/$\Delta t$) & Conflicts Points ($|\rho_v|$)\\
\midrule
Right & Priority & 10\% & 5 & 4 & 5 & 2 \\
Through & Priority & 80\% & 5 & 4 & 5 & 6 \\
Left & Yield & 10\% & 5& 4 & 5 & 6 \\
\bottomrule
\end{tabular}
\caption{Turning proportions, lane capacity, unconditional movement service rates and number of conflict points (blue phase only) based on movement type. The time period is $\Delta t = 10 s$ and green phase lost time is $L = 2 s$.}
\label{tab:mov}
\end{table}

Prior to solving model \green, it is necessary to identify the potential conflict sets $\C_{ij}$ of each movement $(i,j)$ of the intersection. In this illustration, we assume the following intersection configuration:
\begin{itemize}
\item \textbf{Right movements} are conflicting with through and left movements ending at the same lane.
\item \textbf{Through movements} are conflicting with right movements ending on the same lane, perpendicular through movements and left movements starting from a different lane.
\item \textbf{Left movements} are conflicting with right movements ending on the same lane, through movements starting from different lanes and left movements starting from perpendicular lanes. 
\end{itemize}

To illustrate this intersection configuration, the conflict sets $\C_{ij}$ of movements emanating from lane $S_-$ are provided below (the conflict sets of movements emanating from other lanes are symmetrical):
\begin{itemize}
\item $\C_{S_-,E_+} = \left\{(W_-,E_+), (N_-,E_+)\right\}$
\item $\C_{S_-,N_+} = \left\{(E_-,N_+), (W_-,E_+), (E_-,W_+), (E_-,S_+), (N_-,E_+), (W_-,N_+)\right\}$
\item $\C_{S_-,W_+} = \left\{(N_-,W_+), (W_-,E_+), (N_-,S_+), (E_-,W_+), (E_-,S_+), (W_-,N_+)\right\}$
\end{itemize}

\begin{table}%
\centering
\begin{tabular}{lll}
\toprule
Lane ($i$) & Demand ($x_i(t)$)	& AV route choice (blue phase only)\\
\midrule
$S_-$ & 10 & $[N_+,E_+,N_+,N_+,W_+,N_+,N_+,N_+,N_+,N_+]$ \\
$W_-$ & 4 & $[E_+,E_+,S_+,W_+]$\\
$N_-$ & 2 & $[S_+,S_+]$\\
$E_-$ & 7 & $[W_+,W_+,W_+,N_+,S_+,W_+,W_+]$\\
\bottomrule
\end{tabular}

\caption{Lane-queues length and vehicle route choice (for blue phase only).}
\label{tab:demand}
\end{table}

For comparison purposes, we assume identical traffic demand for both LV and AV-lanes. We assume that there is no downstream traffic, hence pressure weights are equal to lane-demand, \emph{i.e.} $\xit = \wit$ for all incoming lanes into the intersection. Lane-demand and AV route choice information is provided in Table \ref{tab:demand}. 

We explore the behavior of \green and \blue in this base case configuration and in the case where movement capacity is doubled, i.e. $C_i = 10$ vehicles per time period (hereby referred to as ``2$\times$Capacity''). The details of the optimal solutions of \green and \blue are summarized in Tables \ref{tab:lane} (lane-based variables) and \ref{tab:movsol} (movement-based variables). In the base case configuration, the optimal objective function value of \green is \textcolor{blue}{$\zgreen = 50$} and that of \blue is $\zblue = 55$. In the 2$\times$Capacity case, the objective of \green is $\zgreen = 104$ and that of blue is $\zblue = 115$.

\begin{table}%
\centering
\scalebox{0.85}{
\begin{tabular}{l lll lll lll lll}
\toprule
& \multicolumn{6}{l}{Base case ($C_i = 5$ veh/$\Delta t$)} & \multicolumn{6}{l}{2$\times$Capacity ($C_i = 10$ veh/$\Delta t$)} \\
\cmidrule(l){2-7} \cmidrule(l){8-13} 
& \multicolumn{3}{l}{\green} & \multicolumn{3}{l}{\blue} & \multicolumn{3}{l}{\green} & \multicolumn{3}{l}{\blue} \\
\cmidrule(l){2-4} \cmidrule(l){5-7} \cmidrule(l){8-10} \cmidrule(l){11-13}
Lane & $\pit$ & $\yit$ & $\wit \yit$ & $\pit$ & $\yit$ & $\wit \yit$ & $\pit$ & $\yit$ & $\wit \yit$ & $\pit$ & $\yit$ & $\wit \yit$ \\
\midrule
$S_-$ & 0.5 & 5.0 & 50.0 & 0.2  & 2.0 & 20.0 & 1.0 & 10.0 & 100.0 & 0.6  & 6.0  & 60.0 \\
$W_-$ & 0.0 & 0.0 & 0.0  & 0.75 & 3.0 & 12.0 & 0.0 & 0.0  & 0.0   & 1.0  & 4.0  & 16.0 \\
$N_-$ & 0.0 & 0.0 & 0.0  & 0.5  & 1.0 & 2.0  & 1.0 & 2.0  & 4.0   & 1.0  & 2.0  & 4.0 \\
$E_-$ & 0.0 & 0.0 & 0.0  & 0.43 & 3.0 & 21.0 & 0.0 & 0.0  & 0.0   & 0.71 & 5.0  & 35.0\\
\bottomrule
\end{tabular}
}
\caption{Optimal lane-based variables obtained by solving \green and \blue for the illustrative instances. For green phase, lane service rates are calculated as $\yit = \xit \pit$. For blue phase, lane service rates are calculated as $\yit = \sum_{v \in \Vit} z_v$ and FIFO blocking effects as $\pit = \yit/\xit$.}
\label{tab:lane}
\end{table}

\begin{table}%
\centering
\scalebox{0.85}{
\begin{tabular}{l llll lll llll lll}
\toprule
& \multicolumn{7}{l}{Base case ($C_i = 5$ veh/$\Delta t$)} & \multicolumn{7}{l}{2$\times$Capacity ($C_i = 10$ veh/$\Delta t$)} \\
\cmidrule(l){2-8} \cmidrule(l){9-15} 
& \multicolumn{4}{l}{\green} & \multicolumn{3}{l}{\blue} & \multicolumn{4}{l}{\green} & \multicolumn{3}{l}{\blue} \\
\cmidrule(l){2-5} \cmidrule(l){6-8} \cmidrule(l){9-12} \cmidrule(l){13-15}
Movement & $\yijt$ & $\mij$ & $\aijt$ & $\bijt$ & $\yijt$ & $\aijt$ & $\bijt$ & $\yijt$ & $\mij$ & $\aijt$ & $\bijt$ & $\yijt$ & $\aijt$ & $\bijt$ \\
\midrule
$(S_-,E_+)$ & 0.5 & 3.5 & 1.0 & 1 & 1.0 & 0.2 & 1 & 1.0 & 8.0  & 1.0  & 1 & 1.0 & 0.1 & 1 \\
$(S_-,N_+)$ & 4.0 & 0.0 & 1.0 & 1 & 1.0 & 0.2 & 1 & 8.0 & 1.0  & 1.0  & 1 & 4.0 & 0.4 & 1 \\
$(S_-,W_+)$ & 0.5 & 3.5 & 1.0 & 1 & 0.0 & 0.0 & 0 & 1.0 & 8.0  & 0.82 & 1 & 1.0 & 0.1 & 1 \\
$(W_-,S_+)$ & 0.0 & 0.0 & 0.0 & 0 & 1.0 & 0.2 & 1 & 0.0 & 0.0  & 0.0  & 0 & 1.0 & 0.1 & 1 \\
$(W_-,E_+)$ & 0.0 & 0.0 & 0.0 & 0 & 2.0 & 0.4 & 1 & 0.0 & 0.0  & 0.0  & 0 & 2.0 & 0.2 & 1 \\
$(W_-,N_+)$ & 0.0 & 0.0 & 0.0 & 0 & 0.0 & 0.0 & 0 & 0.0 & 0.0  & 0.0  & 0 & 1.0 & 0.1 & 0 \\
$(N_-,E_+)$ & 0.0 & 0.0 & 0.0 & 0 & 0.0 & 0.0 & 0 & 0.2 & 8.8  & 0.11 & 1 & 0.0 & 0.0 & 0 \\
$(N_-,S_+)$ & 0.0 & 0.0 & 0.0 & 0 & 1.0 & 0.2 & 1 & 1.6 & 7.4  & 1.0  & 1 & 2.0 & 0.2 & 1 \\
$(N_-,W_+)$ & 0.0 & 0.0 & 0.0 & 0 & 0.0 & 0.0 & 0 & 0.2 & 8.8  & 1.0  & 1 & 0.0 & 0.0 & 0 \\
$(E_-,N_+)$ & 0.0 & 0.0 & 0.0 & 0 & 0.0 & 0.0 & 0 & 0.0 & 0.0  & 0.0  & 0 & 1.0 & 0.1 & 1 \\
$(E_-,W_+)$ & 0.0 & 0.0 & 0.0 & 0 & 3.0 & 0.6 & 1 & 0.0 & 0.0  & 0.0  & 0 & 3.0 & 0.3 & 1 \\
$(E_-,S_+)$ & 0.0 & 0.0 & 0.0 & 0 & 0.0 & 0.0 & 0 & 0.0 & 0.0  & 0.0  & 0 & 1.0 & 0.1 & 1 \\
\bottomrule
\end{tabular}
}
\caption{Optimal value of movement-based variables obtained by solving \green and \blue for the illustrative instances. For green phase, movement service rates are calculated as $\yijt = \pij\xit\pit$ and $\mij$ is calculated by \eqref{eq:slack}. For blue phase, movement service rates are calculated as $\yijt = \sum_{v \in \V_{ij}^n(t)} z_v$.}
\label{tab:movsol}
\end{table}

Optimal lane service rates and pressure-weighted lane service rates are reported in Table \ref{tab:lane}. These results highlight the difference between \green and \blue with the former servicing only 1 lane in the base case and 2 lanes in the 2$\times$Capacity case. In contrast, \blue is capable to service some traffic on all 4 lanes of the intersection thanks to the conflict-point formulation allowing multiple conflicting movements to be activated simultaneously.

In the base case, only the 3 movements from South lanes ($S_-$) are activated by \green. Blocking effects due to FIFO conditions yield $\phi_{S_-}(t)$=0.5 as summarized in Table \ref{tab:lane}. This indicates that only 50.0\% of the demand on lane $S_-$ is able to access the intersection. The service rate of movements is detailed in Table \ref{tab:movsol}. Of the activated movements, all movements have full service capacity. In this example, the high demand of the South lane leads the intersection manager to only activate South-based movements. This is reflected in the movements slacks: the through movement $(S_-,N_+)$ has null slack ($\mu_{S_-N_+} = 0.0$), thus blocking other movements which conflict with this movements. This activation pattern is typical of saturated intersections. In the 2$\times$Capacity case, there is enough capacity to service all South lane demand $\mu_{S_-N_+} = 1.0$, which is sufficient to allow North lane demand to travel through the intersection. In this case, there are no blocking effects on North and South lanes and yield movements admit partial service levels of $\alpha_{S_-W_+}=0.82$ and $\alpha_{N_-E_+}=0.11$.

In contrast to green phase, blue phase activates at least one movement from each incoming lane to the intersection with a total of 6 movements in the base case, and 9 movements in the 2$\times$Capacity case, as indicated in Table \ref{tab:movsol}. All activated movements have fractional service rates ($\aijt$) comprised between 0.2 and 0.6, which is narrower range compared to services rates of activated movements in green phase. Comparing the number of vehicles moved, we find that in the base case \green moves a total of 5 vehicles while \blue moves 9 vehicles. In the 2$\times$Capacity case, \green is able to move 12 vehicles and \blue moves 17 vehicles.

This highlights that at low capacity, the green phase model is substantially constrained. Further, \green can be particularly sensitive to turning ratios $\pij$. For example, if the turning ratios for left turns are null and the corresponding demand shifted to right turns instead, solving the same example in the base case yields a an activation matrix wherein all through and right turn movements from South and North lanes are activated since all of these movements do not conflict with one another. In this case, $\phi_{S_-}=0.5$ and $\phi_{N_-}=1$ since North lane demand can be fully serviced.

The MILPs \green and \blue can be solved to optimality using traditional branch-and-cut-and-bound techniques in near real-time, as shown in Section \ref{num}. Their solution provides the basis for the proposed online network traffic control policy which is introduced in the next section, along with its proof of stability.

%%%%%%%%%%%%%%%%%%%%%%%%%%%%%%%%%%%%%%%%%%%%%%%%%%%%%%%%%%%%%%%%%%%%%%%%%%%%
\section{Hybrid Network Control Policy and Stability}\label{policy}
%%%%%%%%%%%%%%%%%%%%%%%%%%%%%%%%%%%%%%%%%%%%%%%%%%%%%%%%%%%%%%%%%%%%%%%%%%%%

In this section, we present a new network traffic control policy for intersection control combining green (LV-lane restricted) and blue phases (AV-lane restricted) and prove that it maximizes throughput. The proposed network traffic control policy works by repeatedly solving \green and \blue at each time period $t$ based on the network state $\xt$ and combining local (intersection-level), optimal activation matrices into a network-wide activation matrix $\at$. 

%%%%%%%%%%%%%%%%%%%%%%%%%%%%%%%%%%%%%%%%%%%%%%%%%%%%%%%%%%%%%%%%%%%%%%%%%%%%
\subsection{Stability Region}
%%%%%%%%%%%%%%%%%%%%%%%%%%%%%%%%%%%%%%%%%%%%%%%%%%%%%%%%%%%%%%%%%%%%%%%%%%%%

Let $\AM$ be the set of service matrices. A \textit{policy} is a function $\pi: \X \rightarrow \AM$ that chooses a service matrix $\at \in \AM$ for every state $\xt$. We use the concept of \emph{strong stability} to characterize the proposed stochastic queuing process \citep{leonardi2001bounds,wongpiromsarn2012distributed,zaidi2016back}:
\begin{defi}
Let $\bar{t}$ be a time period index. A stochastic queue evolution process is \textit{strongly stable} under policy $\pi$ if and only if there exists $K < \infty$ such that the network state $\xt$ verifies:
\begin{equation}
\limsup_{\bar{t} \rightarrow \infty} \mathbb{E}\left[\frac{1}{\bar{t}}\sum_{t=1}^{\bar{t}} |\xt| \right] < K	
\label{eq:stability_sq}
\end{equation}
\label{defi:stability_sq}
\end{defi}

For brevity, we hereby referred to \emph{strong stability} as \emph{stability}.\\

To show that the proposed policy is stabilizing and define the stability region of the system, we introduce artificial decision variables in \green and \blue MILPs. Let $\A_{ro} = \A_r \cup \A_o$ and let $\git \in [0,1]$ be a decision variable representing the lane service supply on lane $i \in \A_{ro}$. Recall that $\usi = \sum_{j \in \A : (i,j) \in \M} \usij$ is the unconditional lane supply. Consider the lane supply constraints of the form:
\begin{subequations}
\begin{align}
\yit \leq \git \usi &\quad \forall i \in \A_{ro} \label{eq:side1a} \\
\git \leq \bijt  &\quad \forall (i,j) \in \M \label{eq:side1b}
\end{align}
\label{eq:side1}
\end{subequations}

Constraints \eqref{eq:side1} require that $\git$ be at least $\yit/\usi$ and impose that $\git$ be null if any movement emanating from lane $i$ is inactive. Observe that since $\bijt = 0 \Rightarrow \yijt = 0$, FIFO conditions on lanes imply that if there exists a lane $j$ such that $\pij > 0$ and $\bijt = 0$, then $\yit = 0$, which is consistent with Constraints \eqref{eq:side1}. Observe that Constraints \eqref{eq:side1a} are equivalent to the class-specific lane supply constraints:
\begin{subequations}
\begin{align}
\xit \pit \leq \git \usi &\quad \forall i \in \A_l \label{eq:sideG} \\
\sum_{v \in \V_i^n(t)} z_v \leq \git \usi &\quad \forall i \in \A_a \label{eq:sideB}
\end{align}
\label{eq:side2}
\end{subequations}

Lane supply constraints \eqref{eq:side1} can be incorporated in \green by adding \eqref{eq:sideG} and \eqref{eq:side1b} to \eqref{mod:green}. Similarly, \eqref{eq:side1} can be incorporated in \blue by adding \eqref{eq:sideB} and \eqref{eq:side1b} to \eqref{mod:blue} via Equations \eqref{eq:aijblue}, \eqref{eq:bijblue} and \eqref{eq:phiblue} which link variables $\aijt$, $\bijt$ and $\pit$ to \blue, respectively. Observe that variables $\git \in [0,1]$ do not influence the optimal solution \green or \blue MILPs since $\yit \leq \usi$ for all $i \in \A_{ro}$. We henceforth refer to the MILPs with lane supply constraints \eqref{eq:side2} and variables $\git$, $\forall i \in \A_{ro}$ as the extended \green and \blue MILPs. 
 
Let $\Omega_\green$ and $\Omega_\blue$ denote the feasible regions of the extended \green and \blue MILPs, respectively. Further, let $\Omega = \Omega_\green \cup \Omega_\blue$ be the feasible region verifying the constraint sets of both \green and \blue, and let $\conv(\Omega)$ be the convex hull of $\Omega$. We say that an infinite sequence of controls $\{\gt\}$ for $t = 1, 2, \ldots$ is admissible if $\{\gt \in \Omega\}$ for all time periods $t$.

\begin{defi}
\label{admi}
For any lane $i \in \A_{ro}$, let $\ghi = \lim\limits_{\bar{t}\rightarrow\infty}\frac{1}{\bar{t}}\sum_{t=1}^{\bar{t}} \git$ be the average lane supply of $i$ for sequence $\{\gt \in \Omega\}$. An admissible control sequence $\{\gt \in \Omega\}$ accommodates flow $\bm{f}$ if
\begin{equation}
\fii < \ghi \usi \quad \forall i \in \A_{ro}
\label{eq:admi}
\end{equation}
\end{defi}

Definition \ref{admi} states that a flow $\bm{f}$ can be accommodated if, on average, $\ghi \usi$ vehicles can be served on lane $i \in \A_{ro}$. From Proposition 2 of \citet{varaiya2013max}, $\gh \in \conv(\Omega)$ if and only if $\{\gt \in \Omega\}$ is an admissible control sequence. We define the stability region of the system accordingly.

\begin{defi}
Let $\D$ be the set of demand rate vectors such that there exists an admissible control sequence $\{\gt \in \Omega\}$ with $\gh \in \conv(\Omega)$ defined as:
\begin{equation}
\D \equiv \left\{\bm{d} \in \Re_+^{|\A_r|} : (\fii \leq \ghi \usi, i\in\A_{ro}, \left\{\gt \in \Omega \right\}) \right\}
\label{eq:stabregion}
\end{equation}

Let $\R$ be the interior of $\D$ \emph{i.e.}: $\R = \{\bm{d} \in \D^\circ\}$. $\R$ is the stability region of the system.
\end{defi}

%%%%%%%%%%%%%%%%%%%%%%%%%%%%%%%%%%%%%%%%%%%%%%%%%%%%%%%%%%%%%%%%%%%%%%%%%%%%
\subsection{Hybrid Max-pressure Network Control Policy}
%%%%%%%%%%%%%%%%%%%%%%%%%%%%%%%%%%%%%%%%%%%%%%%%%%%%%%%%%%%%%%%%%%%%%%%%%%%%

The proposed hybrid max-pressure policy selects at each time period $t$ the service matrix $\at$ which maximizes the network-wide pressure among all possible green and blue phases.

\begin{defi}
The network traffic control policy $\pi^\star (\xt)$, defined as:
\begin{equation}
\pi^\star (\xt) = \argmax \left\{\sum_{n \in \N} \big(\zgreen \vee \zblue\big) : \at \in \AM\right\}
\label{eq:opt}
\end{equation}
is hereby referred to as the hybrid max-pressure network control policy for legacy and autonomous vehicles.
\end{defi}

At each intersection $n \in \N$, the hybrid max-pressure policy selects the phase (green or blue) and service matrix ($\at$) maximizing the local pressure. The service matrix is either $\atng$ if the green phase is activated or $\atnb$ if the blue phase is activated. The implementation of the hybrid max-pressure policy requires the resolution of the two MILPs \green and \blue at each intersection of the network. We next show that policy $\pi^\star$ is stabilizing for any demand rate in the stability region of the system $\R$.

\begin{theorem}
\label{theo}
The stochastic queue evolution process \eqref{eq:queue} is stable under the hybrid max-pressure policy $\pi^\star$ \eqref{eq:opt} for any demand rates vector $\bm{d} \in \R$.
\begin{proof}

To show that the stochastic queue evolution process \eqref{eq:queue} verifies the stability condition \eqref{eq:stability_sq}, we first observe that this process is a discrete-time markov chain (DTMC) and we will show that it satisfies the conditions of Theorem 2 in \citet{leonardi2001bounds}. Specifically, we need only to show that there exists $\epsilon > 0$ such that the drift of a chosen Lyapunov function of the DTMC is upper-bounded. Let $\xt \mapsto \left\|\xt\right\|^2  = \sum_{i \in \A} (\xit)^2$ be the chosen Lyapunov function (equivalent to $V(\cdot)$ in the notation of \citet{leonardi2001bounds}). We show that the Lyapunov drift of the DTMC $\mathbb{E}\left[\left\|\xtt\right\|^2 - \left\|\xt\right\|^2 |\ \xt\right]$ is upper bounded by $-\epsilon |\xt|$ where $|\xt| = \sum_{i \in \A} |\xit|$ under the hybrid max-pressure policy $\pi^\star$ \eqref{eq:opt}.\\

Let $\dxit = \xitt - \xit$. We have: $\left\|\xtt\right\|^2 - \left\|\xt\right\|^2$ = $\left\|\xt + \dxt\right\|^2 - \left\|\xt\right\|^2$ = $2\xt^\intercal \cdot \dxt + \left\|\dxt\right\|^2$. Recall that $\A_{ro} = \A_r \cup \A_o$:
% We next show that $\xt^\intercal \cdot \dxt$ and $\left\|\dxt\right\|^2$ are both upper-bounded.\\
\begin{align*}
\xt^\intercal \cdot \dxt =& - \sum_{j \in \A_{ro}} \xjt \Yjt + \sum_{j \in \A_o} \sum_{i \in \A_{ro}} \xjt \pijt \Yit + \sum_{j \in \A_r} \xjt \djt \\
=& - \sum_{i \in \A_{ro}} \xit \Yit + \sum_{j \in \A_o} \sum_{i \in \A_{ro}} \xjt \pijt \Yit + \sum_{j \in \A_r} \xjt \djt \\
=&   \sum_{i \in \A_{ro}} \Yit \left(-\xit + \sum_{j \in \A_o}\xjt \pijt \right) + \sum_{j \in \A_r} \xjt \djt 		
\end{align*}

Computing the conditional expected value $\mathbb{E}[\xt^\intercal \cdot \dxt |\ \xt]$:
\begin{equation}
\mathbb{E}\left[\xt^\intercal \cdot \dxt |\ \xt\right] = - \sum_{i \in \A_{ro}} \mathbb{E}\left[\Yit |\ \xt\right] \wit + \sum_{j \in \A_r} \xjt \dj 
\label{eq:exp}	
\end{equation}

Further, since $\Yit \geq 0$ and $\pijt \leq 1$, we have:
\begin{equation}
\dxjt = - \Yjt + \sum_{i \in \A} \pijt \Yit + \djt \leq \begin{cases}
\sum_{i \in \A_l : (i,j) \in \M_l} \Yit + \djt \text{ if } j \in \A_l \\
\sum_{i \in \A_a : (i,j) \in \M_a} \Yit + \djt \text{ if } j \in \A_a
\end{cases}
\label{eq:l2}
\end{equation}
		
Let $\bar{S}_i$ be the maximum value of the random variable $\Yit$. If $i \in \A_l$, then $\bar{S}_i$ can be determined based on the unconditional lane service capacities $\usi$; otherwise if $i \in \A_a$, then this bound can be derived from AVs' maximum speed under conflict-free traffic conditions. Let $\bar{S} = \max\{\bar{S}_i : i \in \A\}$ and let $K_j$ be the number of lanes permitted to reach lane $j$, \emph{i.e.} $K_j = |i \in \A_l : (i,j) \in \M_l|$ if $j \in \A_l$ or $K_j = |i \in \A_a : (i,j) \in \M_a|$ if $j \in \A_a$. In addition, let $\bar{D}_i$ be the maximum value of the random variable $\dit$. From \eqref{eq:l2} we get $\dxjt \leq K_j \bar{S} + \bar{D}_j$, which gives the following bound:
\begin{equation}
\mathbb{E}\left[\left\|\dxt\right\|^2 |\ \xt\right] \leq \sum_{j \in \A} (K_j \bar{S} + \bar{D}_j)^2 = K'
\label{eq:Z2}
\end{equation} 

By definition $\mathbb{E}\left[\Yit |\ \xt\right] = \yit$. Hence, combining Equations \eqref{eq:exp} and \eqref{eq:Z2}, the Lyapunov drift of the DTMC can be upper-bounded as:
\begin{align}
\mathbb{E}\left[\left\|\xtt\right\|^2 - \left\|\xt\right\|^2 |\ \xt\right] &= \mathbb{E}\left[2\xt^\intercal \cdot \dxt + \left\|\dxt\right\|^2 |\ \xt\right]\\
&\leq K'- 2\left(\sum_{i \in \A_{ro}} \yit \wit + \sum_{j \in \A_r} \xjt \dj\right) \label{eq:driftbound}
\end{align}

Since by assumption $\bm{d} \in \R$ and the stability region is non-empty, there exists an admissible control sequence $\{\gt \in \Omega\}$ with $\gh \in \conv(\Omega)$, which stabilizes $\bm{d}$. Further, by Theorem 2 of \citet{leonardi2001bounds}, there exists $\epsilon^\dagger > 0$, such that the Lyapunov drift of the DTMC is upper-bounded as follows:
\begin{equation}
\mathbb{E}\left[\left\|\xtt\right\|^2 - \left\|\xt\right\|^2 |\ \xt\right] \leq K' - \epsilon^\dagger \sum_{i \in \A_{ro}} |\xit| = K' - 2 \left(\frac{\epsilon^\dagger}{2} \sum_{i \in \A_{ro}} |\xit| \right)
\end{equation}

We next show that the hybrid max-pressure policy $\pi^\star$ achieves a better Lyapunov value than the admissible control sequence $\{\gt \in \Omega\}$. Specifically, if the hybrid max-pressure policy $\pi^\star$ is used, then the term $\sum_{i \in \A_{ro}} \yit \wit$ in \eqref{eq:driftbound} is maximized through the \green and \blue MILPs. Since $\sum_{j \in \A_r} \xjt \dj$ is a constant with regards to the control variables of the policy, and since $\epsilon^\dagger$ is a bound on the difference of the achieved outflow and the achieved inflow under the admissible control sequence $\{\gt \in \Omega\}$, it must hold that
\begin{equation}
\sum_{i \in \A_{ro}} \yits \wit + \sum_{j \in \A_r} \xjt \dj \geq \sum_{i \in \A_{ro}} \yit \wit + \sum_{j \in \A_r} \xjt \dj \geq \frac{\epsilon^\dagger}{2} \sum_{i \in \A_{ro}} |\xit|
\end{equation}

Hence, using policy $\pi^\star$, the upper bound on the Lyapunov drift of the DTMC given in \eqref{eq:driftbound} becomes:
\begin{align}
\mathbb{E}\left[\left\|\xtt\right\|^2 - \left\|\xt\right\|^2 |\ \xt\right] &\leq K'- 2\left(\sum_{i \in \A_{ro}} \yits \wit + \sum_{j \in \A_r} \xjt \dj\right) \label{eq:driftbound2} \\
&\leq K' - \epsilon^\dagger \sum_{i \in \A_{ro}} |\xit| = K' - \epsilon^\dagger |\xt| \label{eq:driftbound3}
\end{align}

Since $K'$ is a constant, Equation \eqref{eq:driftbound3} is equivalent to Condition (3) in Theorem 2 of \citet{leonardi2001bounds}, with the choice of $\xt \mapsto \left\|\xt\right\|^2$ as the Lyapunov function $V()$ and $B = 0$, which proves the theorem.

\end{proof}
\end{theorem}

A natural extension of Theorem \ref{theo} is that the pure network traffic control policies wherein only \green or \blue is used to coordinate traffic are also stable.

\begin{corollary}
\label{coro}
The pure pressure-based network traffic control policy consisting of policy $\pi^\star$ \eqref{eq:opt} with only green (respectively, blue) phases coordinated by \green (respectively, \blue) are stable for any demand rates in the stability region $\R$.

\begin{proof}
To prove that pure network traffic control policies are stable within the stability region $\R$ it suffices to observe that policy $\pi^\star$ \eqref{eq:opt} is defined based on a logical OR-condition taking the maximum local pressure among green and blue phases at each time period and intersection. Hence, Theorem \ref{theo} also applies to the unary case wherein only \green or \blue is used to coordinate traffic.
\end{proof}
\end{corollary}

Theorem \ref{theo} proves that the proposed hybrid network control policy $\pi^\star$ \eqref{eq:opt} stabilizes any demand vector in the stability region of the system $\R$ \eqref{eq:stabregion}. According to \citet{tassiulas1992stability}, and as discussed in \citet{wongpiromsarn2012distributed} and \citet{varaiya2013max} for the case of signalized traffic control, this is equivalent to throughput optimality since it shows that the stability region of policy \eqref{eq:opt} is a superset of the stability region of any policy.

Corollary \ref{coro} establishes that pure policies based on \green or \blue traffic control models also maximize throughput. We note that stability of the pure green network traffic control case is an extension of the work of \citet{varaiya2013max}. \citet{varaiya2013max} proposed a network traffic control policy for a single class of vehicles and assumed that each movement had a dedicated queue. In addition, it was assumed that movement capacities are exogenous to the traffic signal control policy. We have both relaxed this framework by only requiring knowledge of lane-queues, and extended the formulation to two classes of lanes. Further, we introduced pressure-based formulations for both green and blue phases which account for FIFO blocking effects on upstream lanes, and account for the loss of capacity due to conflicting movements being simultaneously activated.

%%%%%%%%%%%%%%%%%%%%%%%%%%%%%%%%%%%%%%%%%%%%%%%%%%%%%%%%%%%%%%%%%%%%%%%%%%%%
\subsection{Online Network Traffic Control Algorithm}
%%%%%%%%%%%%%%%%%%%%%%%%%%%%%%%%%%%%%%%%%%%%%%%%%%%%%%%%%%%%%%%%%%%%%%%%%%%%

We are now ready to present our decentralized network traffic control algorithm used to implement the proposed hybrid max-pressure network control policy. The pseudo-code of the proposed policy is summarized in Algorithm \ref{algo:policy}. At each time period $t$, we calculate the optimal \green and \blue phases at each intersection of the network $n \in \N$ based on the current state of the network $\xt$. The phase with the highest local pressure is selected for each intersection.

\begin{algorithm}
\KwIn{$\G=(\N,\A)$, $\bm{d}$, $\bm{p}$, $\bar{\bm{s}}$, $t$, $\xt$}
\KwOut{$\at$}
\For{$n \in \N$}{
$\zgreen \gets$ Solve \green \eqref{mod:green} \\
$\zblue \gets$ Solve \blue \eqref{mod:blue} \\
$\atn \gets \argmax_{\atng, \atnb} \{\zgreen, \zblue\}$
}
$\at \gets [\atn]_{n\in\N}$
\caption{Hybrid max-pressure network control policy}
\label{algo:policy}
\end{algorithm}

%%%%%%%%%%%%%%%%%%%%%%%%%%%%%%%%%%%%%%%%%%%%%%%%%%%%%%%%%%%%%%%%%%%%%%%%%%%%
\section{Numerical Experiments}\label{num}
%%%%%%%%%%%%%%%%%%%%%%%%%%%%%%%%%%%%%%%%%%%%%%%%%%%%%%%%%%%%%%%%%%%%%%%%%%%%

In this section, we conduct numerical experiments to test the proposed hybrid network control policy and report our findings. 

%%%%%%%%%%%%%%%%%%%%%%%%%%%%%%%%%%%%%%%%%%%%%%%%%%%%%%%%%%%%%%%%%%%%%%%%%%%%
\subsection{Implementation Framework}
%%%%%%%%%%%%%%%%%%%%%%%%%%%%%%%%%%%%%%%%%%%%%%%%%%%%%%%%%%%%%%%%%%%%%%%%%%%%

We implement the proposed hybrid network control policy on artificial datasets to test computational performance and analyze the algorithm's behavior. We use a synthesized grid network of size $5 \times 5$, wherein each of the 25 nodes corresponds to a controlled intersection and each edge represents a bidirectional link between adjacent nodes. All intersections have the same topology as that depicted in Figure \ref{fig:interall}, \emph{i.e.} each node has four incoming and four outgoing links, each of which has one LV-lane and one AV-lane. Each incoming lane allows three movements: through, left turn and right turn. 

We assume that vehicles' routes in the network are fixed. In each instance generated, we randomly and uniformly assign an origin and a destination to each vehicle, a route among these nodes and a departure time within the considered time horizon. Origins and destinations are chosen among nodes at the edge of the grid. The level of travel demand is determined by the \emph{departure rate} of vehicles into the network and the impact of travel demand onto network performance is assessed through a sensitivity analysis.

The time period is set to $\Delta t = 10$ s. We assume that green phases have a \emph{lost time} of 2 s to account for vehicle start-up delays and signal clearance intervals and we conduct a sensitivity analysis on this input parameter. In turn, blue phases are assumed to have zero lost time. 

We use point-queues for all links in the network and we assume that vehicles take three time periods to travel from an intersection to the next intersection on their route. Vehicles travel through 3 links, each with a 10 s free flow time, between each intersection. Hence, in this configuration, it takes 30 s for a vehicle to travel between two adjacent intersections at free flow.  We set the time horizon to 30 minutes and we execute Algorithm \ref{algo:policy} periodically until all vehicles have exited the network.

Vehicles' speed limit through intersections is assumed to be uniform and equal to 44 ft/s and the wave propagation speed is taken as 11 ft/s. We assume that all vehicles have a length of 17.6 ft and that lanes have a width of 12 ft. We use the triangular fundamental diagram to determine lane capacity which results in 1,440 veh/h or 5 vehicles per time period.

We explore the sensitivity of the proposed policy with regards to the proportion of AVs by varying the proportion of AVs from 0\% to 100\% in increments of 10\%. In all our numerical experiments we assume that AVs always choose AV-lanes. To benchmark the performance of the proposed hybrid network control policy (summarized in Algorithm \ref{algo:policy}), we also simulate network traffic under a traditional traffic signal configuration wherein AV-lanes and blue phases are nonexistent. Under this configuration all AV-lanes are treated as LV-lanes, and we model this by using single-lane links with twice the lane capacity of the LV-lanes in the network configuration with both LV- and AV-lanes. Thus, in this benchmark there are no AV-lanes and since we assume that AVs behave as LVs and LV-lanes, all vehicles are treated as LVs. This benchmark is hereby referred to as \textsf{2$\times$Green}. Each experiment is simulated 40 times and average performance is reported.

The simulation framework is implemented in Java on a Windows machine with 8 Gb of RAM and a CPU at 3 GHz. All MILPs are solved with CPLEX 12.8 (Java API) with a time limit of 60 s and default options.

The impact of the departure rate onto network performance is explored in Section \ref{deprate}, the impact of \emph{lost time} during green phases is assessed in Section \ref{losttime} and the activation pattern of green and blue phases is discussed in Section \ref{phase}.

\subsection{Impact of departure rate}\label{deprate}

\begin{figure}%
\begin{center}
\includegraphics[width=0.8\columnwidth]{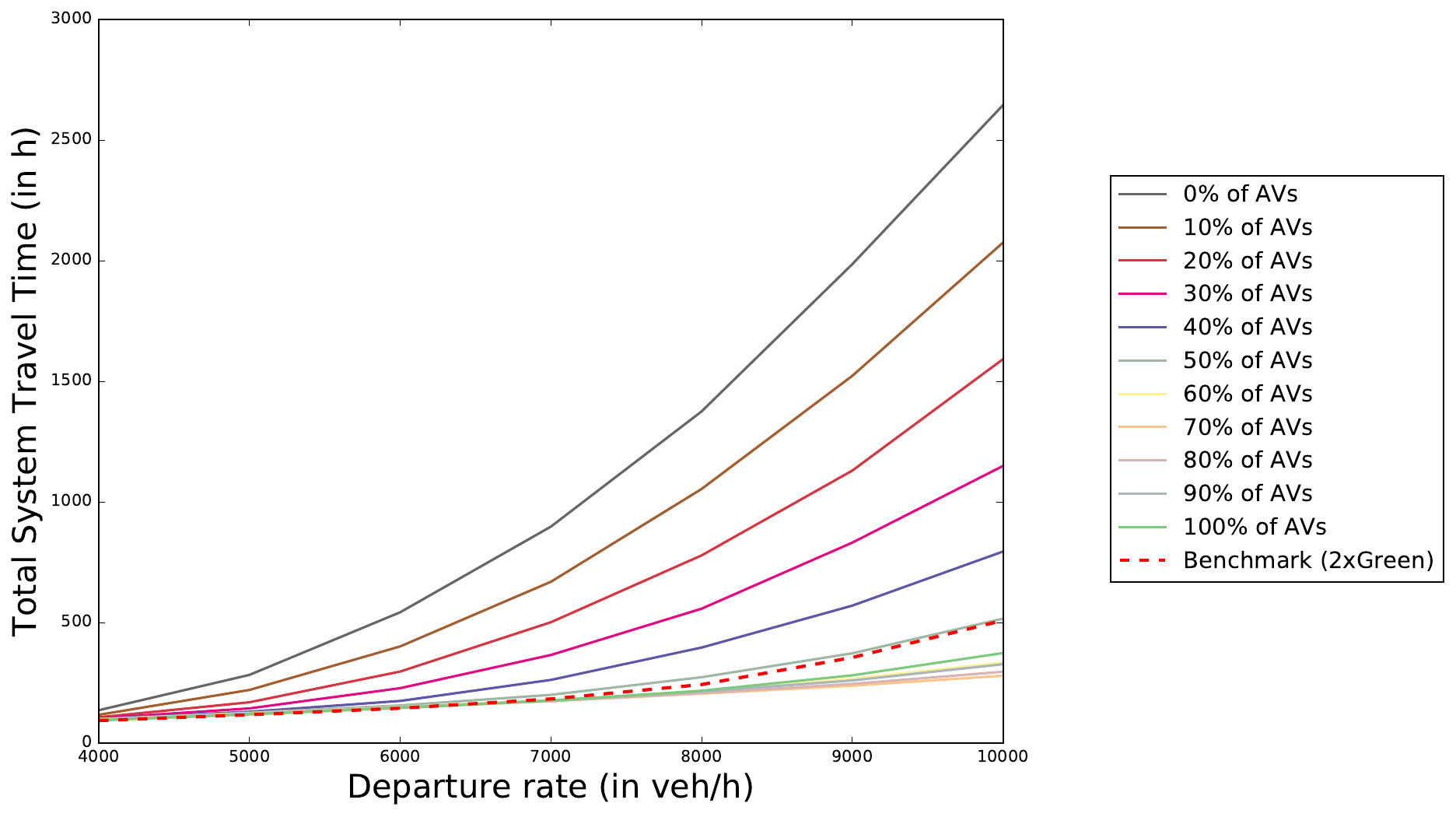}%	
\end{center}
\caption{Total system travel time for a varying proportion of AVs and vehicle departure rate. The results illustrate the trend of the mean total system travel time over 40 simulations on a $5 \times 5$ grid network with each link having one AV and one LV lane. The \textsf{2$\times$Green} experiment (benchmark) corresponds to the scenario where each link has the capacity of two LV-lanes.}%
\label{fig:exp3tstt}%
\end{figure}

\begin{figure}%
\includegraphics[width=\columnwidth]{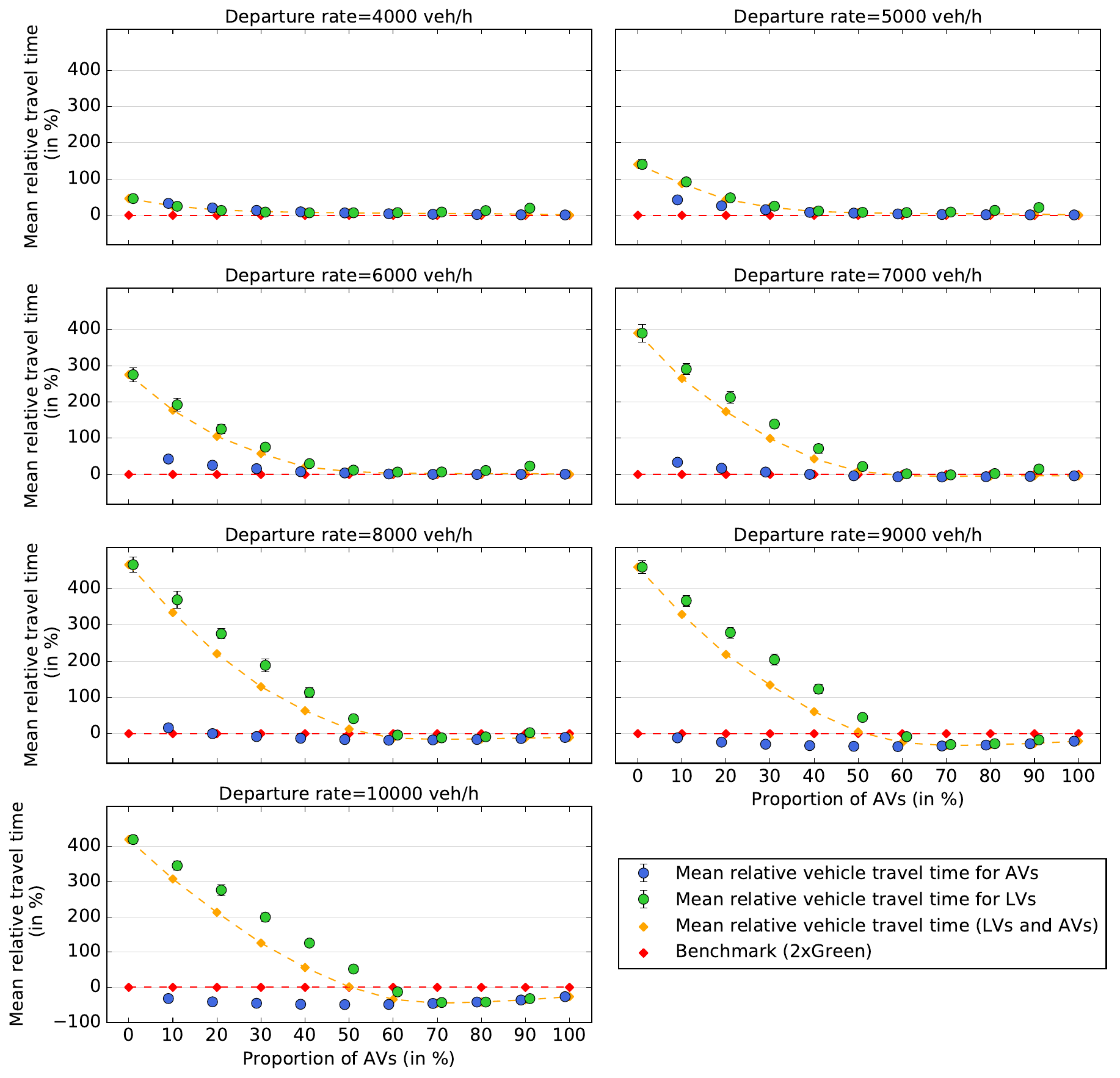}%
\caption{Average vehicle travel time based on AV proportion and vehicle departure rate. The results depict the mean and standard deviation over 40 simulations on a $5 \times 5$ grid network with each link having one AV and one LV lane. The \textsf{2$\times$Green} experiment (benchmark) corresponds to the scenario where each link has the capacity of two LV-lanes.}%
\label{fig:exp3TT}%
\end{figure}

\begin{figure}%
\centering
\includegraphics[width=0.8\columnwidth]{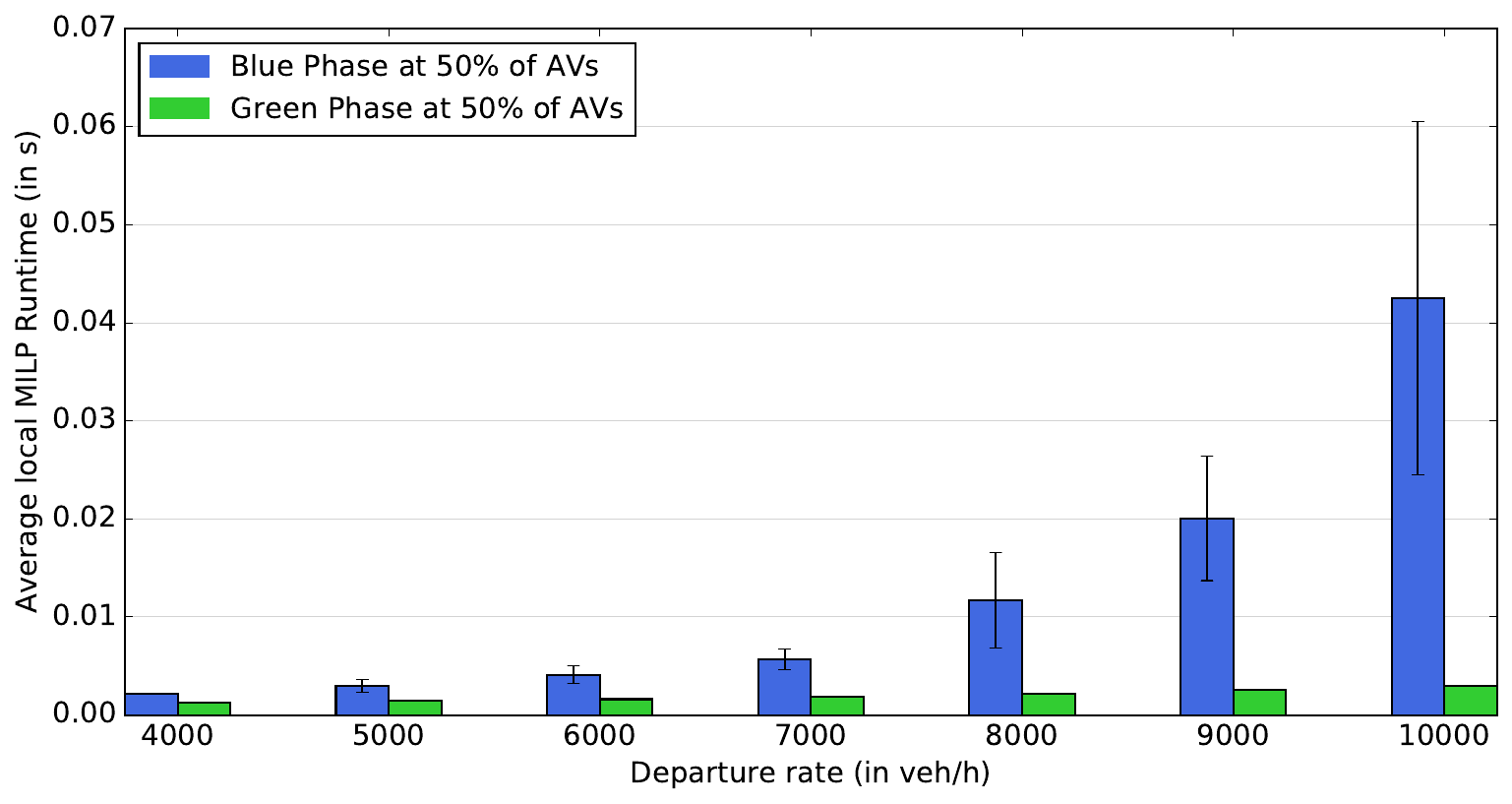}%
\caption{Average, local MILP runtime to solve the \blue or \green to optimality against departure rate (error bars represent standard deviation) over all 25 intersections in the network and time periods required in the simulation. The proportion of AVs is 50\% and the lost time for green phase is 2s.}%
\label{fig:runtime}%
\end{figure}

The evolution of the total system travel time (TSTT) for a varying departure rate is depicted in Figure \ref{fig:exp3tstt}. For this experiment, the green phase lost time is set to 2 s and we vary the departure rate from 4,000 veh/h (lowest demand) to 10,000 (highest demand). As expected, we observe that TSTT increases super-linearly with the departure rate, \emph{i.e.} travel demand. We find that the market penetration of AVs has a significant impact on TSTT. If less than 50\% of AVs are present in the system, we find that the use of dedicated AV lanes and the blue phase is not beneficial for the network in terms of TSTT, even at high demands, as noted by the performance of the benchmark which outperforms the hybrid network control policy for these levels of AV market penetration. Recall that the benchmark represents the TSTT when the pure, green network traffic control policy is used with the equivalent of two LV-lanes. At high demands (\emph{i.e.} 8,000 veh/h and beyond), we find that a market penetration of more than 50\% of AVs outperforms the benchmark. At the highest departure rate tested (\emph{i.e.} 10,000 veh/h), for a market penetration of 70\%, the average TSTT is reduced by 1/2 compared to the benchmark. For AV market penetration rates of 80\% and higher, we find that the reduction of the average TSTT is lesser than that achieved at 70\% of AV market penetration. This can be explained by the fact that network capacity is not fully used at 100\% of AVs compared to 70\% since in the former configuration LV-lanes remain empty throughout the experiment.\\

To further investigate the behavior of the proposed hybrid network control policy, we examine average, vehicle-class travel times relative to the benchmark configuration. Figure \ref{fig:exp3TT} shows the average vehicle travel time for AVs (blue series), LVs (green series) and overall (orange series) in the network based on AV proportion and departure rate. The benchmark is shown as a dashed flat line in red. Three main trends can be identified: first, we find that increasing the departure rate mainly impacts the mean relative travel time of LVs with regards to the benchmark. Second, the proportion of AVs minimizing the mean vehicle travel time in the network is relatively robust to the travel demand and is found to be around 70\% of AVs. This suggests that at lower or higher market penetration rates, network capacity is not saturated and AV- or LV-lanes are under-utilized, respectively. This insight can help in identifying the optimal AV market penetration rate to deploy AV lanes and blue phases. Third, for high levels of travel demand, a sufficiently high proportion of AVs improves on the benchmark, \emph{i.e.} the average vehicle travel time (over both LVs and AVs) obtained using the hybrid network control policy is lower than the average vehicle travel time obtained with the pure, green network control policy.

For a departure rate of 4,000 veh/h, LVs and AVs' average travel time remain similar and the hybrid network control policy performs similarly to the benchmark. Increasing the departure rate from 5,000 to 7,000 veh/h, we observe congestion effects impacting LVs' average travel time, while AVs' average travel time remain only slightly penalized. Further increasing the departure rate to 8,000 veh/h and beyond yields a new pattern: for a proportion of AVs greater or equal to 60\%, the aggregate average vehicle travel time improves on the benchmark although LVs' average travel time remain more penalized than that of AVs. At a departure rate of 10,000 veh/h, we find that both LVs and AVs' average travel time are lower than the average travel time in the benchmark configuration and considerable travel time reduction are achieved from 60\% of AVs and beyond. Specifically, at this departure rate and with a market penetration of AVs of 70\%, we find that the average vehicle travel time is reduced by approximately 1/2 compared to the benchmark. We also observe that the mean relative vehicle (LVs and AVs) travel time exhibits a convex-shaped profile with a minima at 70\% of AV market penetration. This highlights that for high departure rates, high AV market penetration rates (80\% and above) yield an imbalanced usage of network capacity due to the low demand of LVs. \\

The computational performance of the proposed hybrid network control policy is illustrated in Figure \ref{fig:runtime}, which depicts the average runtime of MILPs \blue and \green over all intersections and all simulations against the departure rate. The results are reported for a market penetration corresponding to 50\% of AVs. Computational runtime increases linearly with departure rate for the \green MILP. The \blue MILP computational performance profile exhibits a super-linear growth with departure rate. This is expected since the \blue MILP models each vehicle trajectory whereas the \green MILP use aggregate flow variables to model vehicles movements. Nevertheless, all MILPs are solved in a few milliseconds. Further, the performance MILPs appears to be robust to vehicles' route choice and departure time as demonstrated by the low variance of the computational runtime over intersections and simulations. In practice, since the proposed hybrid network control policy is decentralized, the system is easily scalable to arbitrary-size networks.

\subsection{Impact of green phase lost time}\label{losttime}

\begin{figure}%
\begin{center}
\includegraphics[width=0.8\columnwidth]{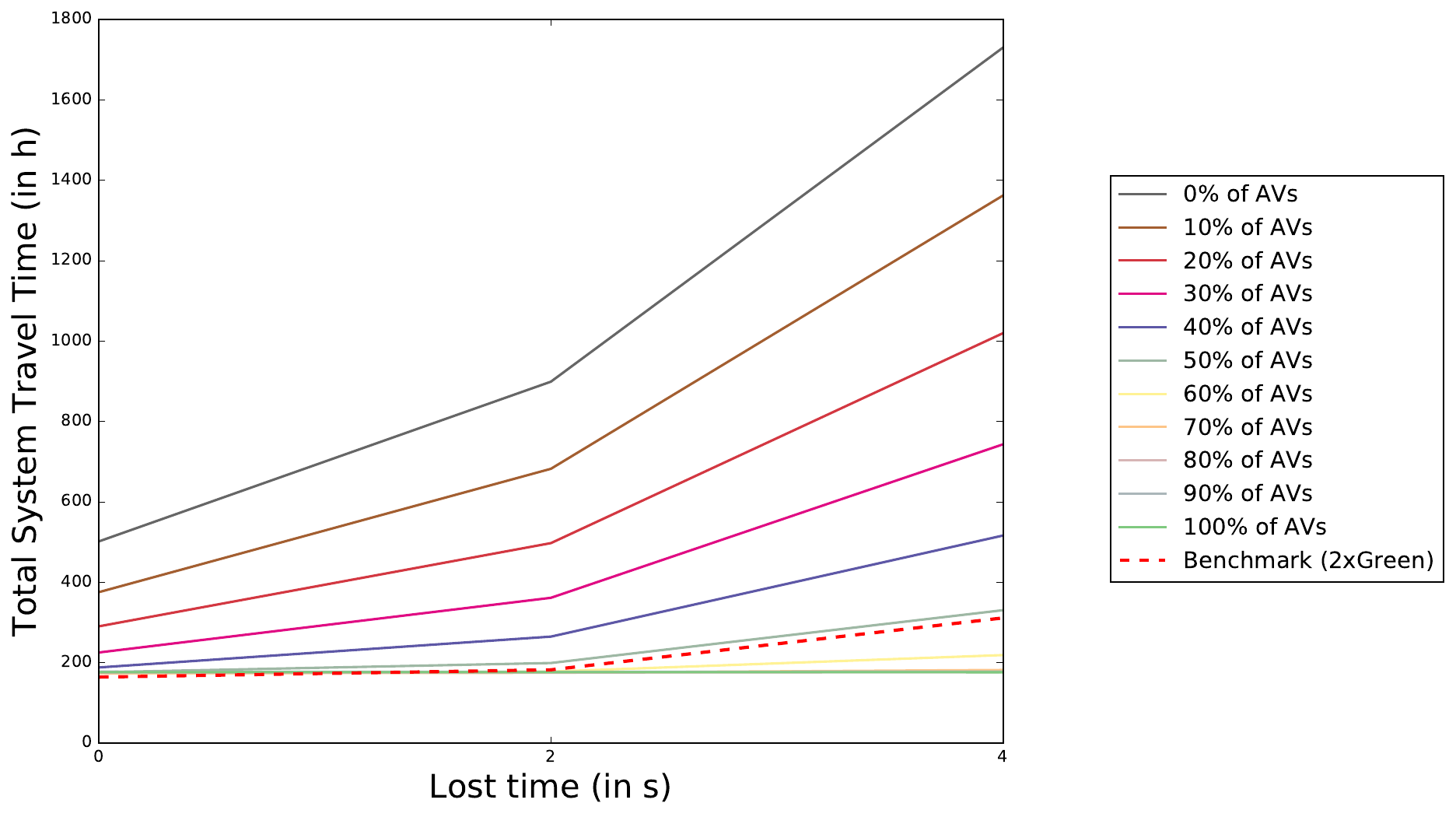}%	
\end{center}
\caption{Total system travel time based on AV proportion and green phase lost time for a departure rate of 7,000 veh/h. The results illustrate the trend of the mean total system travel time over 40 simulations on a $5 \times 5$ grid network with each link having one AV and one LV lane. The \textsf{2$\times$Green} experiment corresponds to the scenario where each link has the capacity of two LV-lanes.}%
\label{fig:exp1tstt}%
\end{figure}

We next explore the impact of green phase lost time. For this sensitivity analysis we set the departure rate to 7,000 veh/h and compare the baseline configuration (lost time of 2 s) with the cases of null (0 s) and doubled (4 s) lost time. The impact of green phase lost time on TSTT is illustrated in Figure \ref{fig:exp1tstt}. We observe that green phase lost time has a super-linear effect on TSTT. For a high level of AV market penetration, \emph{i.e.} with 80\% of AVs or more, the resulting traffic configuration is almost insensitive to lost time. We also find that the hybrid network control policy outperforms the benchmark for a sufficiently high market penetration of AVs (in this case, a proportion 60\% of AVs). 

\begin{figure}%
\includegraphics[width=\columnwidth]{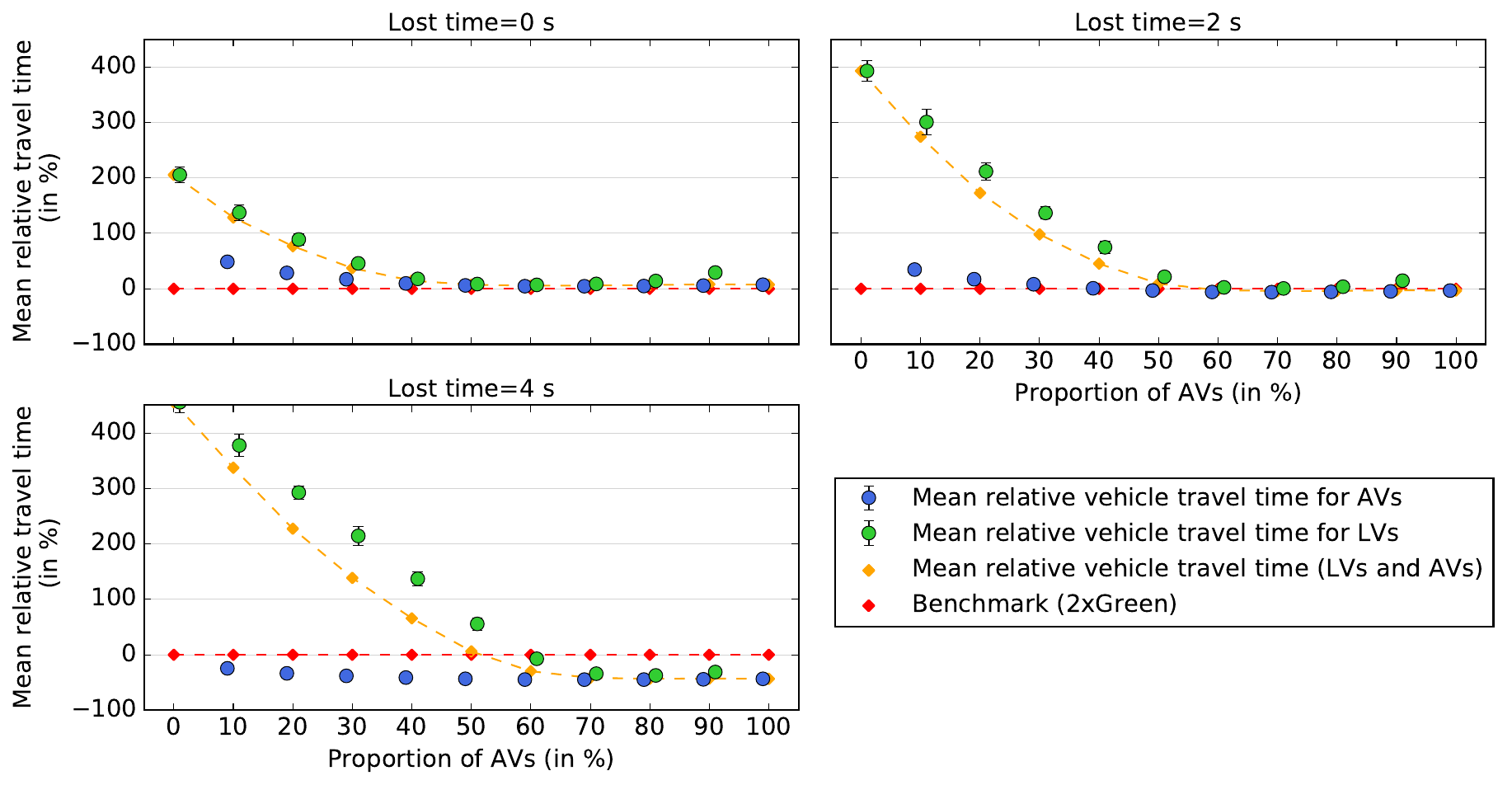}%
\caption{Average vehicle travel time based on AV proportion and green phase lost time for a departure rate of 7,000 veh/h. The results depict the mean and standard deviation over 40 simulations on a $5 \times 5$ grid network with each link having one AV and one LV lane. The \textsf{2$\times$Green} experiment correspond to the scenario where each link has the capacity of two LV-lanes.}%
\label{fig:exp1tt}%
\end{figure}

Figure \ref{fig:exp1tt} provides a more detailed outlook on the impact of green phase lost time by examining vehicle-class average travel time. If the green phase lost time is assumed null, we find that LVs' average travel time remain comparable to AVs' average travel time, especially if both classes of vehicles are in similar proportions in the network, \emph{i.e.} for proportions of AVs between 50\% to 70\%. In turn, doubling the baseline lost time (4 s) considerably penalizes LVs' average travel time for low proportions of AVs. For an AV market penetration of more than 60\%, increasing green phase lost time is better managed with the proposed hybrid network control policy than the in the pure, green network traffic control policy. This can be attributed to the traditional green phase model used in the benchmark configuration: since all lanes are LV-lanes, increasing lost time for green phases results in a considerable loss of capacity compared to the hybrid network configuration.

\subsection{Impact of vehicle spacing at conflict points in blue phase}\label{reserv}

We next explore the impact of Vehicle spacing at conflict points, which can be viewed as a parameter of \blue. Specifically, we incorporate a coefficient $\xi$ in front of $\tau_v(c)$ variables, resulting in the following modification to constraint \eqref{eq:tau}:
\begin{align}
\tau_v(c) = \xi\left(\frac{L_v}{w} + \frac{L_v(t_v(\out_v)-t_v(\inc_v))}{d_v(\inc_v,\out_v)}\right) && \forall v \in \Vt, \forall c\in\p_v
\end{align}

\begin{figure}%
\begin{center}
\includegraphics[width=0.8\columnwidth]{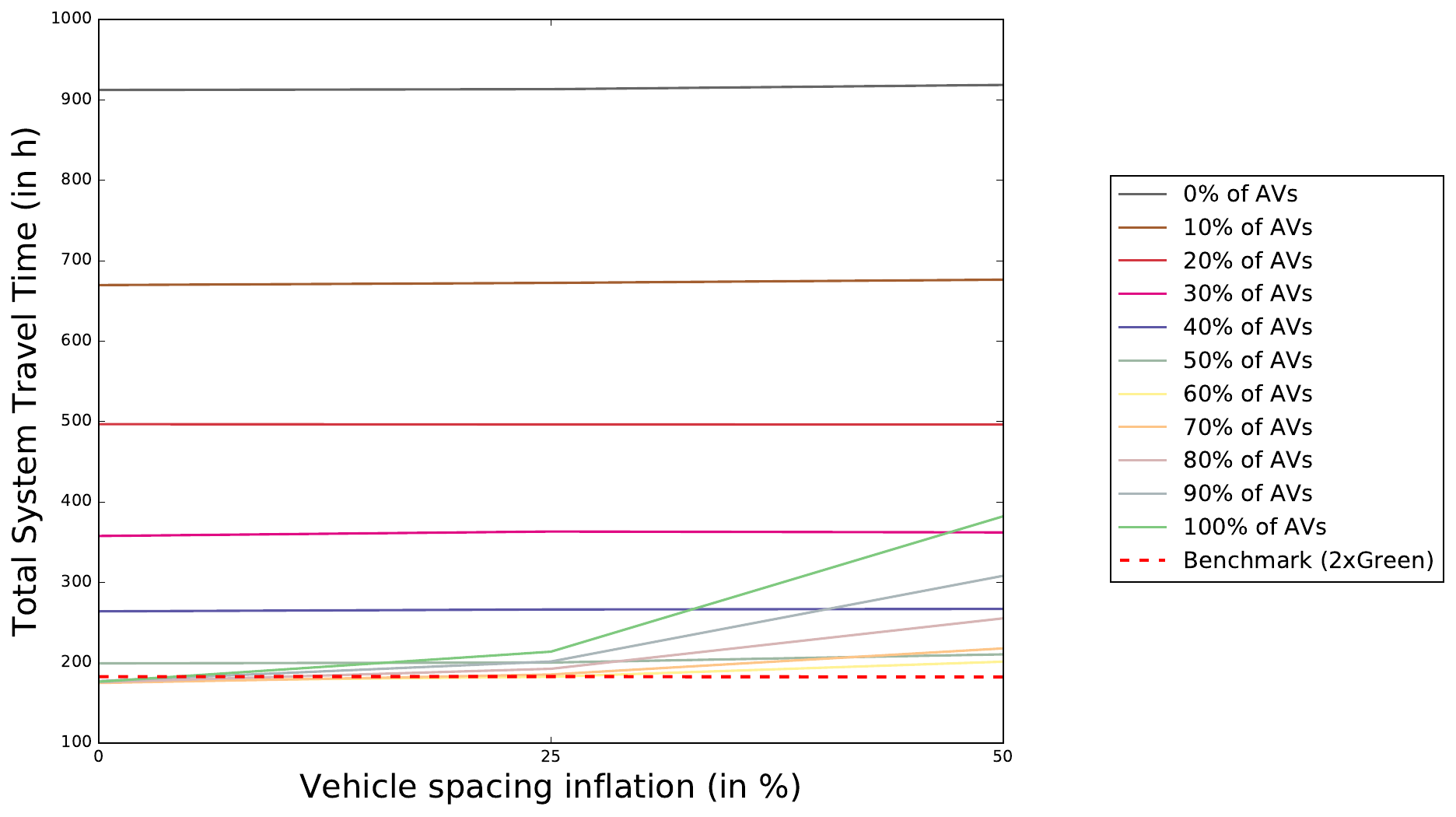}%	
\end{center}
\caption{Total system travel time based on AV proportion and vehicle spacing inflation (blue phase) for a departure rate of 7,000 veh/h. The results illustrate the trend of the mean total system travel time over 40 simulations on a $5 \times 5$ grid network with each link having one AV and one LV lane. The \textsf{2$\times$Green} experiment corresponds to the scenario where each link has the capacity of two LV-lanes.}%
\label{fig:exp4tstt}%
\end{figure}

\begin{figure}%
\includegraphics[width=\columnwidth]{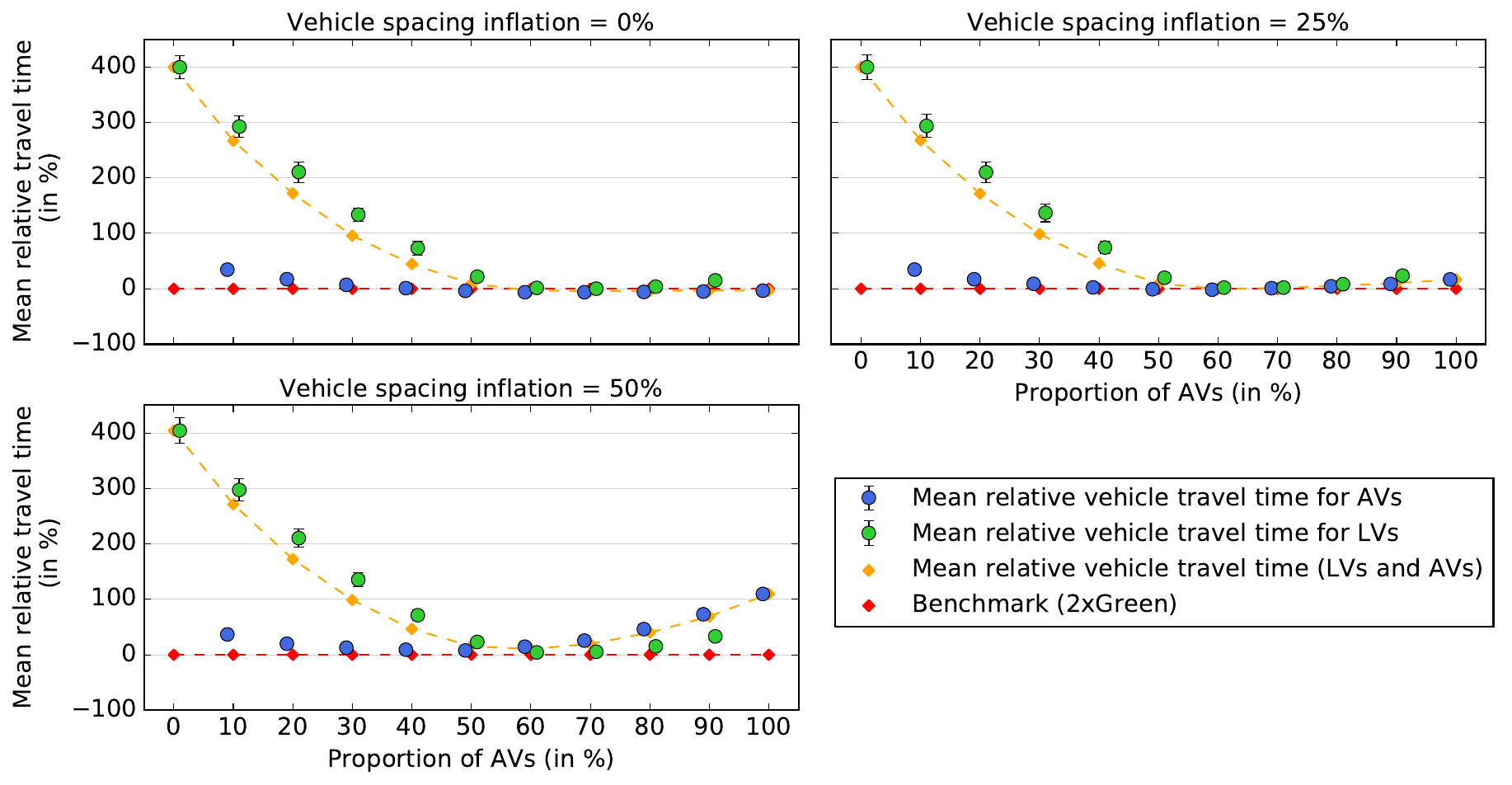}%
\caption{Average vehicle travel time based on AV proportion and vehicle spacing inflation (blue phase) for a departure rate of 7,000 veh/h. The results depict the mean and standard deviation over 40 simulations on a $5 \times 5$ grid network with each link having one AV and one LV lane. The \textsf{2$\times$Green} experiment correspond to the scenario where each link has the capacity of two LV-lanes.}%
\label{fig:exp4tt}%
\end{figure}

We vary $\xi$ from 1 (baseline) to 1.25 (25\% inflation) and 1.5 (50\% inflation). For this sensitivity analysis we set the departure rate to 7,000 veh/h and the green phase lost time to 2 s. The results are reported in Figure \ref{fig:exp4tstt} (TSTT) and \ref{fig:exp4tt} (class-based vehicle travel time). We observe that inflating $\xi$ in \blue increases the TSTT for AV market penetration rates of 50\% and higher. This suggests that lower levels of AVs in the network, blue phases are not restricting system performance. Comparing with the benchmark, we observe that increasing $\xi$ in blue phases results in higher delays, which can be expected since in the baseline case the proposed hybrid policy only matches the benchmark. In terms of vehicle-class impacts, we find that inflating $\xi$ mainly penalizes AVs at high AV market penetration rates. In particular, for an inflation of 50\% and an AV market penetration rate of 70\% or higher, we find that the mean LV travel time is lower than the mean AV travel time, which contrasts to the no-inflation case wherein AVs have better service.\\

Both this analysis and the analysis of the impact of green phase lost time (\ref{losttime}) offer insights on how the parameters of the system can be tuned to increase safety standards, e.g. by increasing lost time in green phases or vehicle spacing in blue phases. Further, these analyses also suggest that the parameters of the system could be used to prioritize a vehicle class over another. For instance, increasing vehicle spacing in blue phases can favor LV travel time. Conversely, increasing lost time in green phase can favor AV travel time. Such prioritization schemes could play key roles in the incorporation of hybrid traffic control policies in urban transport networks.

\subsection{Phase activation patterns}\label{phase}

\begin{figure}%
\includegraphics[width=\columnwidth]{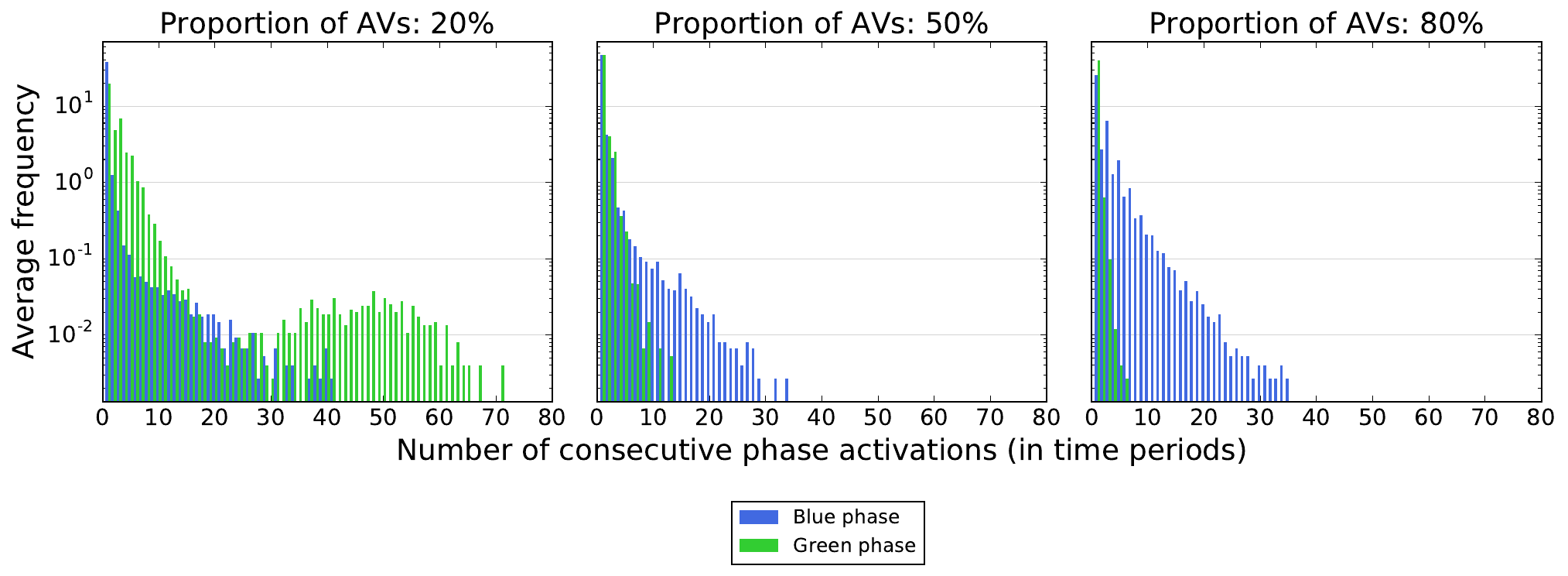}%
\caption{Distribution of average consecutive phase activations over time periods per intersection and per simulation. The departure rate is 7,000 veh/h and the lost time is set to 2 s.}%
\label{fig:phase}%
\end{figure}

To further analyze the behavior of the proposed hybrid network control policy, we report the number of consecutive green or blue phases activated based on the market penetration of AVs in the network. For this analysis, we focus on a departure rate of 7,000 veh/h and a green phase lost time of 2 s, which exhibited the most balanced level of congestion for the generated instances and three AV penetration rates: 20\%, 50\% and 80\%. Figure \ref{fig:phase} depicts the distribution of average consecutive phase activation over time periods per intersection and per simulation. For an AV market penetration of 20\%, we find that green phase is activated for more consecutive time periods compared to blue phase, with up to $\sim$70 consecutive activations. This figure considerably reduces when the proportion of AVs in the network increases. Increasing the proportion of AVs to 50\% and 80\% results in more frequent short streaks for green phases with the majority of consecutive activations between 1 and 10 time periods. In contrast, blue phase activations are found to be more robust to the levels of AVs in the network. Blue phases are mainly activated in streaks of 1 to 30 time periods, with a majority of consecutive activations between 1 and 10 time periods. Since the blue phase admits more combinations of vehicle movements compared to the green phase, the latter often does not have greater pressure until queue lengths are longer. Hence, if the level of AVs in the network is high, LVs may have to wait several time periods before being serviced.

%%%%%%%%%%%%%%%%%%%%%%%%%%%%%%%%%%%%%%%%%%%%%%%%%%%%%%%%%%%%%%%%%%%%%%%%%%%%
\section{Discussion and Perspectives}\label{con}
%%%%%%%%%%%%%%%%%%%%%%%%%%%%%%%%%%%%%%%%%%%%%%%%%%%%%%%%%%%%%%%%%%%%%%%%%%%%

In this paper, we proposed a new, pressure-based network traffic control policy for coordinating legacy vehicles (LV) and autonomous vehicles (AV). The proposed approach assumes that LVs and AVs share the urban network infrastructure. Specifically, we hypothesized that dedicated AV-lanes are available for AVs to access traffic intersections, thus obviating the inherent limitations of first-in-first-out (FIFO) lane queues, \emph{i.e.} vehicle blocking. This design assumption is plausible if the level of penetration of AVs is sufficiently high \citep{levin2016multiclass}. It should be noted that AV-lanes need not to be added infrastructure. If the proportion of AVs in the network is high enough, some LV-lanes can be restricted to AV-traffic. To coordinate traffic at network intersections, we introduced two intersection-level MILP formulations for maximizing local pressure. \green is used to coordinate traffic among LV-lanes. The proposed formulation for green phases only requires knowledge of lane queues and conflict-free movement capacities and estimates actual movement capacity endogenously based on movement activation. In addition, since route choice is assumed unknown for LVs, MILP \green accounts for vehicle-blocking effects due to FIFO conditions on lane-queues. To manage AV-traffic at network intersections, we introduce a so-called \emph{blue} phase during which only AVs are allowed to access the intersection. The resulting \blue MILP is adapted for max-pressure control from the conflict-point formulation introduced by \citet{levin2017conflict}. We characterized the stability region of the proposed queuing system and showed that the proposed decentralized hybrid network control policy is stable, \emph{i.e.} that it maximizes network throughput, under conventional travel demand conditions. Further, Theorem 1 and its proof show how traffic control formulations which are based on vehicle- or trajectory-level variables can be incorporated in stable network traffic control policies.

We conducted numerical experiments on randomly generated artificial instances on a grid-network to test the proposed policy. We explored the sensitivity of the policy with regards to the proportion of AVs in the network as well as the departure rate, which corresponds to the level of travel demand. We also investigated the impact of green phase lost time and examined consecutive phase activation patterns. We found that the different patterns emerged based on the level of congestion in the network. At low congestion levels, we observe that AVs' travel time remains close to the benchmark travel time whereas LVs' travel time is increasingly penalized when the proportion of AVs exceeds that of LVs. A low proportion of AVs penalizes AVs' travel time; instead a high proportion of AVs penalizes LVs. At higher congestion levels, the hybrid network control policy is seen to outperform the benchmark, thus quantifying the benefits the blue phase model for coordinating AV-traffic. We also find that travel demand mainly impacts LVs' travel time whereas AVs' travel time are considerably less penalized. Identifying critical levels of penetration for AVs can provide insight into the management of urban infrastructure. For instance, this can help in assessing at which point it becomes beneficial to restrict specific lanes to AV-traffic.\\

The outcomes of the experiments also reveal that fairness should be taken into consideration when allocating lanes to vehicle classes (LVs, AVs). Indeed, for high proportions of AVs in the network, LVs' travel time may be considerably penalized, despite the overall average travel time improving. Exploring the trade-offs between network throughput and fairness will be addressed in future studies. Real-time lane allocation among vehicle classes (\emph{e.g.} LV, AV) can be expected to further impact route choice and travel behavior altogether. Specifically, while LVs may be assumed to be restricted to LV-lanes for safety and traffic stability reasons; AVs may have the choice between using AV-lanes or LV-lanes. Lane choice behavior may have strategic implications for AVs and thus the proposed network traffic control policy. For instance, if AV-lanes are tolled, AVs with low value of time may choose to avoid them whereas high value of time users may prefer them. This more general network design and control problem can be modeled using bilevel or simulation-based optimization wherein users' departure time and route choice can be accounted for based on traffic equilibrium theory \citep{le2017utility}. Given that our study is focused on traffic control and assumes fixed route choice, we leave this investigation for future research. 

The proposed hybrid network control policy presents some limitations which outline possible future research directions. Notably, we have assumed a point-queue model which does not account for physical link capacities. As discussed by Gregoire et. al. \cite{gregoire2014capacity} physical link capacities can be incorporated in pressure-based policies but may compromise their stability properties. Further, if at a given intersection and time period, all pressure weights are null or negative, then pressure-based control policies have no incentive to move traffic downstream. This is a well-known limitation of such policies which may result in loss of work conservation, i.e. no vehicles being moved at a time period while having nonzero queue lengths. However, The comparison of the max-pressure solution to the average intersection actuation (which by assumption provides sufficient capacity for stability) analytically shows that the max-pressure control indeed achieves a better Lyapunov value. Further, scaling-up pressure weights so that they are always positive is not trivial in point-queue models, which allow infinitely large queues. Hence, further research is needed to establish stabilizing pressure-based control policies which ensure no loss of work conservation and are able to account for physical link capacities.
 
From a practical standpoint, the proposed model admits arbitrary phase ordering and unbounded cycle lengths, which may not be acceptable in practice. While this is a common design limitation in the literature of max-pressure control \citep{wongpiromsarn2012distributed,varaiya2013max,xiao2014pressure}, future work should investigate the effects of these limitations on throughput and maximum stability in max-pressure control. 
 
The treatment of both classes of lanes also presents opportunities for improvement. The proposed hybrid network control policy is designed to treat LV- and AV-lanes equally, i.e. no class is given any form of priority. This design choice could be relaxed to create incentives for network optimization. For instance, in a context where so-called shared-AVs dominate the AV market, promoting AV-lanes could be an attractive strategy to reduce car ownership and transition towards shared mobility. In turn, at low AV penetration rates, it may be beneficial to prioritize LV-traffic to compensate for the reduced capacity of green phase compared to blue phases. This design choice may require significant modification to the stability analysis, necessitating a new max-pressure policy.

%%%%%%%%%%%%%%%%%%%%%%%%%%%%%%%%%%%%%%%%%%%%%%%%%%%%%%%%%%%%%%%%%%%%%%%%%%%%
\bibliography{MWL}
\bibliographystyle{abbrvnat}
%%%%%%%%%%%%%%%%%%%%%%%%%%%%%%%%%%%%%%%%%%%%%%%%%%%%%%%%%%%%%%%%%%%%%%%%%%%%

\end{document}